\documentclass[11pt]{article}
\usepackage{epsfig}
\usepackage{amsthm}
\usepackage{amssymb}
\usepackage{amsmath}
\usepackage{amsfonts}
\usepackage{color}
\newcommand{\field}[1]{\mathbb{#1}}
\newcommand{\R}{\field{R}}

\newcommand{\N}{\field{N}}
\newcommand{\Z}{\field{Z}}

\newcommand{\E}{\field{E}}

\def\qed{\hfill$\diamondsuit$}

\theoremstyle{example} \theoremstyle{remark} \theoremstyle{lemma}
\theoremstyle{definition} \theoremstyle{corol}
\theoremstyle{proposition} \theoremstyle{condition}
\theoremstyle{assumption}
\newtheorem{assumption}{\n{Assumption}}[section]
\newtheorem{theorem}{\n{Theorem}}[section]
\newtheorem{example}{\n{Example}}[section]
\newtheorem{remark}{\n{Remark}}[section]
\newtheorem{lemma}{\n{Lemma}}[section]

\font\n=cmcsc10
\def\cov{{\mbox{cov}}}
\def\var{{\mbox{var}}}
\def\cum{{\mbox{cum}}}
\begin{document}

\bibliographystyle{plain}

\centerline{\Large  TESTING FOR WHITE NOISE UNDER UNKNOWN}

\centerline{\Large  DEPENDENCE AND ITS APPLICATIONS}

\centerline{\Large TO GOODNESS-OF-FIT FOR TIME SERIES MODELS \footnote{I would like to
thank Professor Pentti Saikkonen and two referees for constructive comments that led to improvement
of the paper. The work is supported in part by NSF grant DMS-0804937. Address correspondence to: Xiaofeng Shao, Department of
Statistics, University of Illinois at Urbana-Champaign,  725 South
Wright St, Champaign, IL, 61820; e-mail: xshao@uiuc.edu}}

\bigskip
 \centerline{\textsc{By Xiaofeng Shao}}
\centerline{\today} \centerline {\it University of Illinois at
Urbana-Champaign}

\bigskip

Testing for white noise  has been well studied in the literature of
econometrics and statistics. For most of the proposed test statistics, such as the
well-known Box-Pierce's  test statistic with fixed lag truncation number, the asymptotic null
distributions  are obtained under independent and identically
distributed assumptions and may not be valid for the dependent white
noise.  Due to recent popularity of conditional heteroscedastic models (e.g. GARCH models),
which imply nonlinear dependence with zero autocorrelation, there is a need to understand the asymptotic properties of the
existing test statistics under unknown dependence. In this paper, we
showed that the asymptotic null distribution of Box-Pierce's test statistic with general weights still holds under unknown
 weak dependence so long as  the lag truncation number grows at an appropriate rate with increasing sample size.
  Further applications to diagnostic checking of the ARMA
 and FARIMA models with dependent white noise errors are also addressed.
 Our results go beyond earlier ones by allowing non-Gaussian and conditional heteroscedastic errors in the ARMA and
 FARIMA models and provide theoretical support for some empirical findings reported in the literature.

 %As an important technical contribution, we partially solve
 %an earlconjecture by
\noindent

 \pagenumbering{arabic}

\setcounter{page}{1}
\section{Introduction}

A fundamental problem in time series analysis is to test for white
noise (or lack of serial correlation). For a   zero-mean stationary process
$\{u_t\}$ with finite variance $\sigma^2=\var(u_t)$, denote its
covariance and correlation functions by $R_u(k)=\cov(u_t,u_{t+k})$
and $\rho_u(k)=R_u(k)/\sigma^2, k\in\Z$ respectively.
 Then the null and alternative hypothesis are
 \[H_0:~\rho_u(j)=0~\mbox{for all}~j\not=0,~\mbox{and}~H_1:~\rho_u(j)\not=0~\mbox{for some}~j\not=0.\]
 Let $f_u(\lambda)=(2\pi)^{-1}\sum_{k\in\Z}\rho_u(k)e^{ik\lambda}$
be the normalized spectral density function of $u_t$. The
equivalent frequency domain expressions to $H_0$ and $H_1$ are
\begin{eqnarray*}
\label{eq:spectral} H_0: f_u(w)=\frac{1}{2\pi},~~w\in
[-\pi,\pi)~\mbox{and}~H_1: f_u(w)\not=\frac{1}{2\pi},~\mbox{for
some}~w\in [-\pi,\pi).
\end{eqnarray*}
In statistical modeling, diagnostic checking is an integrable part
of model building. A common way of testing the adequacy of the
proposed model is by checking the assumption of white noise
residuals. Systematic departure from this assumption implies the inadequacy of
the fitted model. Thus testing for white noise is an
important research topic and it has been extensively studied in
 the literature of econometrics  and statistics.

The methodologies can be roughly divided into two categories: time
domain tests and frequency domain tests. In the time domain, the
most popular test is probably  Box and Pierce's (1970) (BP)
portmanteau test, which admits the following form:
\[Q_n=\sum_{j=1}^{m}\hat{\rho}^2_u(h),\]
where $m$ is the so-called lag truncation number [see Hong (1996)] and is (typically) assumed to be fixed. The empirical autocorrelation  $\hat{\rho}_u(j)$,  is defined as
$\hat{\rho}_u(j)=\hat{R}_u(j)/\hat{R}_u(0)$ with
$\hat{R}_u(j)=n^{-1}\sum_{t=|j|+1}^{n}(u_t-\bar{u})(u_{t-|j|}-\bar{u})$,
where $\bar{u}=n^{-1}\sum_{t=1}^{n}u_t$.
 Under the assumption that $\{u_t\}_{t\in\Z}$ are independent and identically distributed (iid),
 it can be shown that $nQ_n\rightarrow_{D}\chi^2(m)$, where
$``\rightarrow_{D}"$ stands for convergence in distribution. If
$\{u_t\}_{t=1}^{n}$ are replaced by the residuals from a well specified model, then
the limiting distribution is still $\chi^2$ but the degree of
freedom is reduced to $m-m'$, where $m'$ is the number of parameters in the
model. In the frequency domain, Bartlett (1955) proposed test statistics based on the famous
$U_p$ and $T_p$ processes and a rigorous theoretical
treatment of their limiting distributions was provided by  Grenander
and Rosenblatt (1957). Other contributions to the frequency domain
tests can be found in Durlauf (1991) and Deo (2000) among others.

In the literature, when deriving the asymptotic null distribution of the test statistic, most earlier works assume Gaussianity and thus
lack of correlation is equivalent to independence. Lately there
has been work that stress the distinction between lack of
correlation and independence. The main reason is that the
asymptotic null distributions of the above-mentioned test
statistics were obtained under iid assumptions on $u_t$, and may not hold in the presence
 of nonlinear dependence, such as
conditional heteroscedasticity. For example, Romano and Thombs
(1996) showed  that the BP statistic with $\chi^2$ approximation can
lead to misleading inferences when the time series is uncorrelated
but dependent. Francq et al. (2005) also demonstrated that the BP
test applied to the residuals of an ARMA model with uncorrelated
but dependent errors performs poorly without suitable
modifications. Various methods have been proposed to account for the dependence;
see for example, Romano and Thombs (1996),  Lobato et al. (2002),
Francq et al. (2005) and Horowitz et al. (2006) among others.
At this point, it seems natural to ask: ``Does there exist a test
statistic whose asymptotic null distribution is robust to the unknown dependence of $u_t$". We
shall  give an affirmative answer in this paper.

In a seminal paper,  Hong (1996) proposed several test statistics,
which measure the distance between a kernel-based spectral density
estimator and the spectral density of the noise under the null
hypothesis. Let
\[\hat{f}_n(w)=(2\pi)^{-1}\sum_{j=-n+1}^{n-1}K(j/m_n)\hat{\rho}_u(j)e^{ijw}\]
be the lag window estimator of the normalized spectral density
function [Priestley (1981)], where $K(\cdot)$ is a nonnegative
symmetric kernel function, $m_n$ is the bandwidth that depends on
the sample size. With the quadratic distance, Hong's statistic is
expressed as
\[T_n=\pi\int_{-\pi}^{\pi}(\hat{f}_n(w)-(2\pi)^{-1})^2dw,
\]
or equivalently,
\begin{eqnarray*}
T_n=\sum_{j=1}^{n}K^2(j/m_n)\hat{\rho}^2_{u}(j).
\end{eqnarray*}
It is worth noting that BP statistic can be regarded as a special
case of Hong's, where $K(\cdot)$ is taken to be the truncated
kernel $K(x)={\bf 1}(|x|\le 1)$. Under the iid assumptions on
$u_t$ and $1/m_n+m_n/n\rightarrow 0$, Hong (1996) established the
asymptotic null distribution of $T_n$, i.e.
\begin{eqnarray}
\label{eq:mainresult}
\frac{nT_n-C_n(K)}{\sqrt{2D_n(K)}}\rightarrow_{D} N(0,1),
\end{eqnarray}
where $C_n(K)=\sum_{j=1}^{n-1}(1-j/n)K^2(j/m_n)$,
$D_n(K)=\sum_{j=1}^{n-2}(1-j/n)(1-(j+1)/n)K^4(j/m_n)$ and $N(0,1)$ stands for the standard normal distribution.
 Under some
additional assumptions on $K(\cdot)$ and $m_n$,
(\ref{eq:mainresult}) holds with $C_n(K)$ and $D_n(K)$ replaced by
$m_nC(K)$ and $m_nD(K)$ respectively, where
$C(K)=\int_{0}^{\infty}K^2(x)dx$ and
$D(K)=\int_{0}^{\infty}K^4(x)dx$. Later Hong and Lee (2003)
established the above result assuming $u_t$ to be martingale
differences with conditional heteroscedasticity of unknown form. One of the
major contributions of this paper is to show that Hong's test
statistic is still asymptotically valid under general white noise
assumption on $u_t$.  Further, we establish that when replacing
$u_t$ by $\hat{u}_t$, the residuals from the ARMA model
with uncorrelated and dependent errors, the asymptotic null distribution of $T_n$
still holds. Our assumptions and results differ from those in Francq et al. (2005) in that $m$ is held fixed
in their asymptotic distributional theory, while $m=m(n)$ grows with the sample size $n$ in our setting.
From a theoretical standpoint, the fourth cumulant of $u_t$ plays a non-negligible role in the asymptotic distribution of $Q_n$ when $m$ is fixed,
whereas it turns out to be asymptotically negligible in $T_n$ when $m_n\rightarrow\infty$. So in the latter case, the asymptotic null distribution does not change under dependent white noise, i.e. the dependence is automatically accounted for
if $m$ and $n$ both grow to infinity. The theoretical finding is also consistent with the empirical results reported in the simulation studies of
Francq et al. (2005), where the empirical size of the BP test is seen to be reasonably close to the nominal one when $n$ is large and $m$ is relatively large compared to  $n$.

 Recently, there has been considerable attention
paid to the goodness-of-fit for long memory time series. Here we
only mention some representative works. Extending Hong's (1996)
idea, Chen and Deo (2004a) proposed a generalized portmanteau test
based on the discrete spectral average estimator and obtained the
asymptotic null distribution for Gaussian long memory time series.
Following the early work by
Bartlett (1955), Delgado et al. (2005)  studied Bartlett's $T_p$
process with estimated parameters and  a martingale transform approach was
 used to make the null distribution asymptotically
distribution-free. In a related work, Hidalgo and Kreiss (2006)
proposed to use bootstrap methods in the frequency domain to approximate
the sampling distribution of Bartlett's $T_p$ statistic with
estimated parameters. In these two articles, the asymptotic
distributional theory heavily relies on the assumption that the
noise processes are  conditionally homoscedastic martingale
differences.

 In the last decade, the FARIMA (fractional autoregressive integrated moving average) models with GARCH errors have
been widely used in the modeling literature [cf. Lien and Tse (1999), Elek
and M\'{a}rkus (2004), Koopman et al. (2007)]. In the modeling stage of a FARIMA-GARCH
model, it is customary to fit a FARIMA model first and then fit a
GARCH model to the residuals. It is crucial to specify the FARIMA
model correctly since the model misspecification of the
conditional mean often leads to the misspecification of the
GARCH model; see Lumsdaine and Ng (1999). Thus diagnostic
checking of FARIMA models with unknown GARCH errors is a very
important issue.  Note that Ling and Li (1997) and Li and Li (2008) have studied the BP type tests for FARIMA-GARCH models assuming a parametric form for the GARCH model. To the best of our knowledge, there
seems no diagnostic checking methodology known or theoretically
justified to work for long memory time series models with nonparametric conditionally heteroscedastic martingale
difference errors. In this article, we shall fill this gap by proving asymptotic validity of
 Hong's test statistic when we replace the unobserved errors by the estimated counterpart from a FARIMA model.

 We now
introduce some notation. For a column vector $x = (x_1, \cdots,
x_q)'\in \R^q$, let $|x| = (\sum_{j=1}^q x_j^2)^{1/2}$. For a random vector $\xi$, write $\xi \in {\cal L}^p$ ($p > 0$) if
$\|\xi\|_p := [\E(|\xi|^p )]^{1/p} < \infty$ and let $\|\cdot\| =
\|\cdot\|_2$. For $\xi \in {\cal L}^1$ define projection operators
${\cal P}_k \xi = \E( \xi | {\cal F}_k) - \E( \xi | {\cal
F}_{k-1})$, $k \in \Z$, where ${\cal F}_k = (\ldots,
\varepsilon_{k-1}, \varepsilon_k)$ with $\{\varepsilon_t\}_{t\in\Z}$ being iid random variables.
Let $C>0$ denote a generic constant which may vary from
line to line;
denote by ${\rightarrow}_{p}$ convergence
 in probability. The symbols
$O_{p}(1)$ and $o_{p}(1)$ signify being bounded in probability and
convergence to zero in probability respectively.
The paper is structured as follows. In Section \ref{sec:obs} we
introduce our assumptions on $u_t$ and establish the asymptotic
 distributions of $T_n$ under the null and alternative hypothesis.
Section~\ref{sec:unobs} discusses the case when $u_t$ are
not directly observable. Here we consider the ARMA
 and FARIMA models with dependent white noise errors
in Section~\ref{subsec:ARMA} and Section~\ref{subsec:FARIMA}
respectively. Section~\ref{sec:con} concludes.  Proofs are gathered
in Section~\ref{sec:tech}.

\section{When $u_t$ is observable}
\label{sec:obs}

Suitable structural assumptions on the process $(u_t)$ are
certainly needed. Throughout, we assume that $(u_t)$ is a mean
zero stationary causal process of the form
\begin{eqnarray}
\label{eq:ut} u_t=F(\cdots,\varepsilon_{t-1},\varepsilon_t),
\end{eqnarray}
where $\varepsilon_t$ are iid random variables, and  $F$ is a measurable function for which $u_t$ is well defined. Further
we assume $u_t$ satisfies the geometric-moment contraction (GMC) condition
[Hsing and Wu (2004), Shao and Wu (2007), Wu and Shao (2004)]. Let $(\varepsilon_k')_{k\in\Z}$ be an iid copy of
$(\varepsilon_k)_{k\in\Z}$; let
$u_n'=F(\cdots,\varepsilon_{-1}',\varepsilon_0',\varepsilon_1,\cdots,\varepsilon_n)$
be a coupled version of $u_n$. We say that $u_n$ is GMC$(\alpha)$,
$\alpha>0$, if there exist $C>0$ and $\rho=\rho(\alpha)\in (0,1)$
such that
\begin{eqnarray}
\label{eq:gmc} \E(|u_n-u_n'|^{\alpha})\le C\rho^n, ~~n\in\N.
\end{eqnarray}
The property (\ref{eq:gmc}) indicates that the process $\{u_n\}$
forgets its past exponentially fast, and it can be verified for
many nonlinear time series models, such as threshold model, bilinear model
and various forms of GARCH models; see Wu and Min (2005) and Shao and Wu
(2007).

Besides conditional heteroscedastic models, which imply uncorrelation due to the  martingale difference structure, there
are a few commonly used models [see Lobato et al. (2002)] that are uncorrelated but are not martingale differences.
We shall show that these models satisfy GMC property under appropriate assumptions.
\begin{example}
\label{ex:1}
{\rm Bilinear model [Granger and Anderson (1978)]:
 \[u_t=\varepsilon_t+b \varepsilon_{t-1}u_{t-2},\]
  where $\varepsilon_t$ are iid $N(0,\sigma_{\varepsilon}^2)$ and $|b|<1$. According to Example 5.3 in
Shao and Wu (2007), $u_t$ is GMC$(\alpha)$, $\alpha\ge 1$ if
\[\E\left|\left(\begin{array}{cc}
0&1\\
b\varepsilon_t&0
\end{array}\right)
\right|_{\alpha}<1,\]
where for a $p \times p$ matrix $A$, $|A|_\alpha = \sup_{z\not=0} |A
z|_\alpha / |z|_\alpha$, $\alpha \ge 1$, is the matrix norm
induced by the vector norm $|z |_{\alpha} = (\sum_{j=1}^p
|z_j|^\alpha)^{1/\alpha}$.}
\end{example}

\begin{example}
\label{ex:2}
{\rm All-Pass ARMA(1,1) model [Breidt et al. (2001)]:
\[u_t=\phi u_{t-1}+\varepsilon_t-\phi^{-1} \varepsilon_{t-1}\]
where $|\phi|<1$ and $\varepsilon_t\sim iid(0,\sigma_{\varepsilon}^2)$. Note that $u_t=\varepsilon_t+\sum_{j=1}^{\infty}(\phi^j-\phi^{j-2})\varepsilon_{t-j}$. Since $|\phi|<1$, $u_t$ is
GMC$(\alpha)$ if $\varepsilon_t\in {\cal L}^{\alpha}$. In view of Theorem 5.2 in Shao and Wu (2007), the all-pass ARMA$(p,p)$ model also satisfies GMC$(\alpha)$ provided that
 $\varepsilon_t\in {\cal L}^{\alpha}$.
}
\end{example}

\begin{example}
\label{ex:3}
{\rm Nonlinear moving average model [Granger and Ter\"asvirta (1993)]:
\[u_t=\beta\varepsilon_{t-1}\varepsilon_{t-2}+\varepsilon_t,\]
where $\varepsilon_t\sim iid(0,\sigma_{\varepsilon}^2)$ and $\beta\in\R$. It is easily seen that $u_t$ is GMC$(\alpha)$ if $\varepsilon_t\in {\cal L}^{\alpha}$.
}
\end{example}

To obtain the asymptotic distribution of $T_n$, the following
assumption is made on the kernel function $K(\cdot)$ and is
satisfied by several commonly-used kernels in spectral analysis,
such as Bartlett, Parzen and Tukey kernels (see Priestley (1981),
p 446-447).
\begin{assumption}
\label{as:kernel}
 Assume the kernel function $K:\R\rightarrow
[-1,1]$ has compact support on $[-1,1]$, is  differentiable except
at a finite number of points and symmetric with $K(0)=1$,
$\max_{x\in [-1,1]}|K(x)|=K_0<\infty$.
\end{assumption}

The assumption that $K(\cdot)$ has compact support can presumably be relaxed at the expense of longer and more technical
proof; see Chen and Deo (2004a). Here we decide to retain it to avoid more technical complications.

\begin{theorem}
\label{th:null} Suppose Assumption~\ref{as:kernel} and
(\ref{eq:gmc}) holds with $\alpha=8$. Assume $\log n=o(m_n)$ and
$m_n=o(n^{1/2})$. Under $H_0$, we have
\begin{eqnarray}
\label{eq:hongresult}
\frac{nT_n-m_nC(K)}{\sqrt{2m_nD(K)}}\rightarrow_{D} N(0,1).
\end{eqnarray}
\end{theorem}

\begin{remark}
{\rm  As pointed out by a referee,
the $8$-th moment condition on $u_t$ is fairly strong and it excludes some interesting GARCH
models, such as the IGARCH model. In addition, the permissible parameter space for the regular
 GARCH$(r,s)$ model is quite small under the $8$-th moment assumption. At this point, we are unable to relax this assumption as it seems necessary
in our technical argument. Nevertheless, the result above suggests that the asymptotic null
distribution of Hong's (1996) statistic is unaffected by
unknown (weak) dependence. From a technical point of view, the asymptotic null
distribution of the BP statistic depends on the fourth cumulants
of $u_t$ since the number of lags $m$ is fixed.
In contrast, for Hong's
statistic, as $m_n\rightarrow\infty$, the fourth cumulant
effect appears to be asymptotically negligible. For a fixed $m$, our result
in Theorem~\ref{th:null} is not applicable.

 }
\end{remark}

The condition on the bandwidth is less restrictive than it looks.
I am not aware of any theoretical results on the optimal bandwidth
choice for  $T_n$ in the hypothesis testing context.  In terms of
estimating the spectral density function, the optimal bandwidth is
$m_n=Cn^{1/5}$ if the kernel (e.g. Parzen kernel) is quadratic
around zero, and $m_n=Cn^{1/3}$ if the kernel (e.g. Bartlett
kernel) is linear around zero. Note that the problem of testing for white noise
 bears some resemblance to testing lack of fit (or specification testing)
in the nonparametric regression context. The latter problem has been well studied in the literature and
 the data-driven bandwidth choice for the smoothing type test has been addressed in
 Horowitz and Spokoiny (2001) and Guerre and Lavergne (2005) among others.

  For the optimal choice of the
kernel function, we refer the reader to Hong (1996) for more
details. The consistency of $T_n$ is stated in the following
theorem.

\begin{theorem}
\label{th:consistent} Suppose Assumption~\ref{as:kernel} and
(\ref{eq:gmc}) holds with $\alpha=8$. Assume $1/m_n+m_n/n\rightarrow
0$. Under $H_1$, we have
\begin{eqnarray*}
\frac{\sqrt{m_n}}{n}\left(\frac{nT_n-m_nC(K)}{\sqrt{2m_nD(K)}}\right)\rightarrow_{p}\frac12
\sum_{j\not=0}\rho_u^2(j)/(2D(K))^{1/2}.
\end{eqnarray*}
\end{theorem}
\noindent Proof of Theorem~\ref{th:consistent}: It follows from the argument
in the proof of Theorem 6 of Hong (1996) by noting that $R_u(j)\le
Cr^{j}$ for some $r\in [0,1)$ and the absolute summability of the
fourth cumulants under GMC$(4)$ [See Wu and Shao (2004), Proposition
2]. We omit the details. \qed

\begin{remark}{\rm
In a related work, Chen and Deo (2006) considered the variance ratio statistic to test for white noise based on the first
 differenced series and proved that when the horizon $k$ satisfies $1/k+k/n=o(1)$, the asymptotic null distribution of
the variance ratio statistic is also robust to conditional heteroscedasticity of unknown form. Their result
 is akin to ours, in that the asymptotic null distribution of the test statistic is nuisance parameter free and
 the horizon $k$ in variance ratio statistic plays a similar role as our bandwidth $m_n$. However, in their conditions (A1)-(A6), the white noise process is assumed to be a sequence of martingale differences with additional
regularity conditions imposed on the higher order moments (up to $8$th); compare Deo (2000). Under our framework,
 the white noise process does not have to be martingale difference under the null. This has some practical implications since there are nonlinear time series models that are uncorrelated but are not martingale differences, as shown in Examples~\ref{ex:1}-\ref{ex:3}. From a technical point of view, the relaxation of the martingale difference assumption, which was imposed in Hong and Lee (2003) and Chen and Deo (2006), is a very nontrivial step and is made feasible with the novel martingale approximation techniques; see Appendix for more discussions. }

\end{remark}

\begin{remark}
\label{re:BPtest}
{\rm
For the BP test statistic, $K(x)={\bf 1}(|x|\le 1)$ and $C(K)=D(K)=1$. Thus
the statement (\ref{eq:hongresult}) reduces to $\{n\sum_{j=1}^{m_n}\hat{\rho}_{{u}}^2(j)-m_n\}/\sqrt{2m_n}\rightarrow_{D} N(0,1)$.
 In the implementation of the BP test, we use the critical values based on $\chi^2(m_n)$ and compare it with the realized value of $n\sum_{j=1}^{m_n}\hat{\rho}_{{u}}^2(j)$, whereas in Hong's test, the critical values are based on the standard normal distribution. Loosely speaking, the two procedures are asymptotically equivalent, since as $m_n\rightarrow\infty$, the central limit theorem implies
  $\chi^2(m_n)\approx N(m_n,2m_n)$. This suggests that the use of BP test is valid in the presence of unknown weak dependence
  when $m_n$ is relatively large compared to $n$.
}
\end{remark}

\section{When $u_t$ is unobservable}
\label{sec:unobs}
 In practice, the errors $\{u_t\}_{t=1,2,\cdots,n}$ are often unobservable as a part of the model, but can be estimated.
   Hong (1996)
studied the residuals from a linear dynamic model that includes
both lagged dependent variables and exogenous variables. In principle,  our results can be extended to the residuals from any parametric time series models with uncorrelated errors, including the setup studied by Hong (1996). Instead of pursuing full generality,  we shall treat the residuals from ARMA and FARIMA models  in Sections~\ref{subsec:ARMA} and ~\ref{subsec:FARIMA} respectively. This is motivated by the  recent interests on the ARMA models with dependent white
noise errors [cf. Francq and Zako\"1an (2005),
Francq et al. (2005) and the references therein] and goodness-of-fit for long memory time series models [see Section~\ref{subsec:FARIMA} for more references].

\subsection{ {\rm ARMA} model}
\label{subsec:ARMA}
 Consider  a stationary autoregressive and
moving average (ARMA) time series generated by
\begin{eqnarray}
\label{eq:ARMA}
(1-\alpha_1B-\cdots-\alpha_pB^p)X_t=(1+\beta_1B+\cdots+\beta_qB^q)u_t,
\end{eqnarray}
where $B$ is the backward shift operator, $\{u_t\}$ is a sequence
of uncorrelated random variables and
$\Lambda=(\alpha_1,\cdots,\alpha_p,\beta_1,\cdots,\beta_q)'$ is an
unknown parameter vector. Let $\phi_{\Lambda}(z)=1-\alpha_1 z-\cdots-\alpha_p z^p$ and
$\psi_{\Lambda}(z)=1+\beta_1 z+\cdots+\beta_q z^q$ be AR and MA polynomials respectively.
 Denote by
$\Lambda_0=(\alpha_{10},\cdots,\alpha_{p0},\beta_{10},\cdots,\beta_{q0})'$
the true value of $\Lambda$ and assume that $\Lambda_0$ is an
interior point of the set
\[\Omega_{\delta}=\{\Lambda\in\R^{p+q};~\mbox{the roots of
polynomials}~\phi_{\Lambda}(z)~\mbox{and}~\psi_{\Lambda}(z)~\mbox{have
moduli}~\ge 1+\delta\}\] for some $\delta>0$.
Following Francq et al. (2005), we call (\ref{eq:ARMA})
a weak ARMA model if $(u_t)$ is only uncorrelated,  a semi-strong
ARMA model if $(u_t)$ is a martingale difference, and a strong
ARMA model if $(u_t)$ is an iid sequence.

Denote by
$\hat{\Lambda}_n=(\hat{\alpha}_{1n},\cdots,\hat{\alpha}_{pn},\hat{\beta}_{1n},\cdots,\hat{\beta}_{qn})'$
the estimator of $\Lambda$. Then the residuals $\hat{u}_t$,
$t=1,2,\cdots,n$ are usually obtained by the following recursion
\[\hat{u}_{t}=X_t-\hat{\alpha}_{1n}X_{t-1}-\cdots-\hat{\alpha}_{pn}X_{t-p}-\hat{\beta}_{1n}\hat{u}_{t-1}-\cdots-\hat{\beta}_{qn}\hat{u}_{t-q},~t=1,2,\cdots,n,\]
where the initial values
$(X_0,X_{-1},\cdots,X_{1-p})'=(\hat{u}_0,\cdots,\hat{u}_{1-q})'=0$.
Following Francq et al. (2005), we test
\[H_0: (X_t) ~\mbox{has an ARMA}(p,q) ~\mbox{representation}~(\ref{eq:ARMA})\]
against the alternative
\[H_1: (X_t) ~\mbox{does not admit an ARMA representation, or admits an ARMA}(p',q')\]
\vspace{-0.8cm}
\[\mbox{representation with}~ p'>p ~\mbox{or}~ q'>q.\]
If $p$ and $q$ are correctly specified, we would expect the
estimated residuals behave like a white noise sequence under
$H_0$.  The following theorem states the asymptotic null
distribution of the test statistic $T_{1n}=\sum_{j=1}^{n}K^2(j/m_n)\hat{\rho}_{\hat{u}}^2(j)$.

\begin{theorem}
\label{th:armaresidual} Suppose the assumptions in
Theorem~\ref{th:null} hold. Assume
$\hat{\Lambda}_n-\Lambda_0=O_p(n^{-1/2})$. Then under $H_0$,
\[\frac{nT_{1n}-m_nC(K)}{(2m_nD(K))^{1/2}}\rightarrow_{D}N(0,1).\]
\end{theorem}
The proof of Theorem~\ref{th:armaresidual} follows the argument
used in the proof of Theorem~\ref{th:farimaresidual} below and is
simpler. We omit the details. Note that as a common feature of smoothing-type test, the
 use of the residuals $\{\hat{u}_t\}$ in place of the
 true unobservable errors $\{u_t\}$ has no impact on the limiting
 distribution.

\begin{remark}{\rm
In the simulation studies of Francq et al. (2005), it can be seen
that when $m$ is large relative to $n$, the level of the BP test
is reasonably close to the nominal one. Here our result provides  theoretical
support for this phenomenon since if we let $K$ to be the truncated kernel,  the resulting test statistic is
exactly the same as BP's.
As commented in Remark~\ref{re:BPtest}, the difference between the use of the $\chi^2$-based critical values as done in BP test, and
 the use of the $N(0,1)$-based critical values for Hong's test is asymptotically negligible since the number of model parameters (i.e. $p+q$) is fixed and $m_n\rightarrow\infty$. Therefore, it seems fair to say that the use of
 BP test is still justified when the lag truncation number $m$ is large, as the unknown dependence in $u_t$ does not kick in asymptotically.

 %are adopted and the estimation effect is negligible.
%the BP test when $m$ is large, the use of BP test is asymptotically valid [see Remark] and
%the estimation effect is asymptotically negligible.

 }

\end{remark}

As mentioned in Francq et al. (2005), weak ARMA models can arise
from various situations, such as transformation of strong ARMA
processes, causal representation of noncausal ARMA processes and
nonlinear processes. In the sequel, we demonstrate
 that the GMC condition for the noise process in the weak ARMA
respresentation can be verified for the two leading examples in
Francq et al. (2005).
\begin{example}
{\rm Consider the process
\[X_t-aX_{t-1}=\varepsilon_t-b\varepsilon_{t-1},~a\not=b\in (-1,1),\]
where $\varepsilon_t$ are iid random variables with
$\E(\varepsilon_t)=0$ and $\varepsilon_t\in {\cal L}^{\alpha},
\alpha\ge 1$. Let $Y_t=X_{2t}$. Then
$Y_t-a^2Y_{t-1}=\xi_t=u_t-\theta u_{t-1}$, where $\theta\in
(-1,1)$,
$\xi_t=\varepsilon_{2t}+(a-b)\varepsilon_{2t-1}-ab\varepsilon_{2t-2}$,
$u_t$ is  white noise and $u_t=R_{1t}+R_{2t}+\theta
\xi_{t-1}$, where
$R_{1t}=-ab\varepsilon_{2t-2}+\theta^2\varepsilon_{2t-4}+\varepsilon_{2t}+(a-b)\varepsilon_{2t-1}+\theta^2[(a-b)\varepsilon_{2t-5}-ab\varepsilon_{2t-6}]$
and $R_{2t}=\sum_{i\ge 3}\theta^i u_{t-i}$. It is easily seen that
$\xi_t$ and $R_{1t}$ satisfy GMC$(\alpha)$. By Theorem 5.2 in Shao
and Wu (2007), $R_{2t}$ also satisfies GMC$(\alpha)$. Therefore,
$u_t$ is GMC$(\alpha)$.}
\end{example}

\begin{example}{\rm
Consider the process
\[X_t=\varepsilon_t-\phi\varepsilon_{t-1}, ~|\phi|>1.\]
Let $u_t=\sum_{i=0}^{\infty}\phi^{-i}X_{t-i}$. Then $X_t$
admits the causal MA(1) representation:
$X_t=u_t-\phi^{-1}u_{t-1}$. Since $X_t$ is GMC$(\alpha)$, $u_t$ is
also GMC$(\alpha)$ by Theorem 5.2 in Shao and Wu (2007).}
\end{example}

\begin{remark}
{\rm
To study the local power of $T_{1n}$, we follow Hong (1996) and
define the local alternative  $H_{1n}: f_{un}(w)=(2\pi)^{-1}+a_n g(w)$ for $w\in [-\pi,\pi]$, where $a_n=o(1)$.
The function  $g$ is symmetric,  $2\pi$-periodic  and  satisfies
$\int_{-\pi}^{\pi} g(w)dw=0$, which ensures that $f_{un}$ is a valid normalized spectral density function
for large $n$. Let $\mu(K)=2\pi \int_{-\pi}^{\pi} g^2(w) dw/(2D(K))^{1/2}$. It can be shown that
 under $H_{an}$ with $a_n=m_n^{1/4}/n^{1/2}$,
\begin{eqnarray}
\label{eq:la}
\frac{nT_{1n}-m_nC(K)}{(2m_nD(K))^{1/2}}\rightarrow_{D}N(\mu(K),1)
\end{eqnarray}
provided that $\hat{\Lambda}_n-\Lambda_0=O_p(n^{-1/2})$ and the assumptions in Theorem~\ref{th:null} hold.
Since the proof basically repeats the argument in the proof of Hong's (1996) Theorem 4, we omit the details.
It is worth mentioning that the above asymptotic distribution (\ref{eq:la}) under the local alternative still holds for $T_n$, whereas a similar result for $T_{2n}$ [see Section~\ref{subsec:FARIMA} for the definition] in the long memory case may still hold but  the proof seems tedious and is thus not pursued. Compared to the Box-Pierce test with a fixed $m$, Hong's test is locally less powerful in that Box-Pierce's test has nontrivial power against
the local alternative of order $n^{-1/2}$. On the other hand, Box-Pierce's test only has trivial power against non-zero correlations at lags beyond $m$, whereas Hong's test is able to detect non-zero correlations at any nonzero lags asymptotically.
}
\end{remark}

\subsection{{\rm FARIMA} model}
\label{subsec:FARIMA} In this subsection, we extend our result to the
goodness-of-fit problem for long memory time series. A commonly
used model in the long memory time series literature is the FARIMA
model:
 \begin{eqnarray}
 \label{eq:longmemory}
(1-B)^d \phi_{\Lambda}(B)Y_t=\psi_{\Lambda}(B)u_t,
\end{eqnarray}
where $d\in (0,1/2)$ is the long memory parameter.  Let
$\theta=(d,\Lambda')'$ and denote by $\theta_0=(d_0,\Lambda_0')'$
its true value. Assume that $\theta_0$ lies in the interior of
$\Theta_{\delta}=[\Delta_1,\Delta_2]\times \Omega_{\delta}$, where
$0<\Delta_1<\Delta_2<1/2$.

 Testing goodness
of fit for short/long memory time series models has attracted a lot of attention
recently. Most tests were constructed in the frequency domain and they can be roughly categorized into two types:
spectral density based test and spectral distribution function
based test. Tests developed by Hong (1996), Paparoditis (2000),
Chen and Deo (2004a) are of the first type and they usually involve
a smoothing parameter and have trivial power against $n^{-1/2}$
local alternatives. The advantage of this type of tests is that the
asymptotic null distributions are free of nuisance parameters.  For
the second type, see Beran (1992),
Chen and Romano (1999),  Delgado et al. (2005) and Hidalgo and Kreiss
(2006), among others. Typically, the tests of this type avoid the
 issue of choosing the smoothing parameter and they can
distinguish the alternatives within $n^{-1/2}$-neighborhoods of the
null model. However, a disadvantage associated with this kind of tests is that
the asymptotic null distributions often depend on the underlying
data generating mechanism and are not asymptotically distribution-free. The martingale
transform method [see Delgado et al. (2005)] and the bootstrap approach
[Chen and Romano (1999), Hidalgo and Kreiss (2006)] have been
 utilized  to  make the tests practically usable. So far, the tests
proposed by Chen and Deo (2004a), Delgado et al. (2005) and Hidalgo
and Kreiss (2006) have been justified to work for long memory time
series models. However, they  assumed either Gaussian processes or
linear processes with the noise processes being conditionally
homoscedastic martingale differences, which exclude interesting models, such as FARIMA
models with unknown GARCH errors.

  Since $d_0\in (0,1/2)$, the process $Y_t$ is invertible. We have the following autoregressive
 representation
 \[u_t=\sum_{k=0}^{\infty}e_{k}(\theta_0) Y_{t-k}.\]
 Given the observations $Y_t, t=1,2,\cdots,n$, we follow Beran (1995) and
 form the residuals by
 \begin{eqnarray}
 \label{eq:residualfarima}
 \hat{u}_t=\sum_{j=0}^{t-1}e_{j}(\hat{\theta}_n)Y_{t-j},~t=1,2,\cdots,n,
 \end{eqnarray}
 where $\hat{\theta}_n$ is an estimator of $\theta$.
 Similar to the ARMA case, the null and alternative hypothesis are
 \[H_0: (Y_t) ~\mbox{has an FARIMA}(p,d,q) ~\mbox{representation}\]
and
\[H_1: (Y_t) ~\mbox{does not admit an FARIMA representation, or admits an FARIMA}(p',d,q')\]
\vspace{-0.8cm}
\[\mbox{representation with}~ p'>p ~\mbox{or}~ q'>q.\]
The test statistic is
$T_{2n}=\sum_{j=1}^{n}K^2(j/m_n)\hat{\rho}_{\hat{u}}^2(j)$, where $\{\hat{u}_t\}_{t=1}^n$ are from (\ref{eq:residualfarima}).

\begin{theorem}
\label{th:farimaresidual}
 Suppose that the assumptions in Theorem~\ref{th:null} hold.
 Assume $\hat{\theta}_n-\theta_0=O_p(n^{-1/2})$. Then under $H_0$, we have
\[\frac{nT_{2n}-m_nC(K)}{(2m_nD(K))^{1/2}}\rightarrow_{D}N(0,1).\]
\end{theorem}

The result presented above is a new contribution to the literature, even for the model
(\ref{eq:longmemory}) with iid errors. Here  we can take  the  Whittle pseudo-maximum likelihood
estimator as $\hat{\theta}_n$. The root-$n$ asymptotic normality
of Whittle estimator for long memory time series models with general white noise errors
has been established  by Hosoya (1997) and Shao (2010).

\begin{remark}{\rm
Hong's (1996) statistic has been reformulated in the discrete form
by Chen and Deo (2004a), who showed asymptotic equivalence of the two statistics
 for Gaussian long memory time series.
Note that the applicability of Chen and Deo's (2004a) test statistic has only been
proved for the Gaussian case. The latter authors
conjectured that their assumptions can be relaxed to allow long
memory linear processes with iid innovations. The work presented here partially
solves their conjecture and our results even allow for dependent
innovations.

A limitation of our theory is that we need to assume the mean of $Y_t$ is known.
  In practice, if the mean is unknown, we need to modify our $\hat{u}_t$ [cf. (\ref{eq:residualfarima})] by replacing $Y_t$
  with $Y_t-\bar{Y}_n$, where $\bar{Y}_n=n^{-1}\sum_{t=1}^{n} Y_t$. It turns out that
   our technical arguments are no longer valid
   with this modification except for the case $d_0\in (0,1/4)$ with additional restrictions on $m_n$. The main reason is that the sample mean of
   a long memory time series converges to the population mean relatively slowly at the rate of $n^{(1/2-d_0)}$. The larger $d_0$ is, the
     slower it becomes. When $d_0\in [1/4,1/2)$, the effect of mean adjustment becomes asymptotically non-negligible. As pointed out by a referee,
      Chen and Deo's (2004a) frequency domain test statistic  is mean invariant, so no mean adjustment is needed. It might be possible to extend the theory presented in Chen and Deo (2004a) directly to the case of dependent innovations, but such an extension seems very challenging and is beyond the scope of this paper. In the short memory case, i.e. $d_0=0$, the mean adjustment does not affect the asymptotic null distribution of the test statistic $T_{1n}$. In other words,  Theorem~\ref{th:armaresidual} still holds if we use the mean adjusted residuals in the calculation of $T_{1n}$.

      % will not be further considered in this paper.
 %    Consequently, the theoretical justification provided in the proof
 % needs to be considerably changed. I am unable to prove such a result in case of the unknown mean except for $d_0\in (0,1/4)$ and imposing
 % restrictive conditions on $m_n$ (e.g. $m_n/n^{1-4d_0}=o(1)$). Compared to the frequency domain version of Hong's test statistic,
 %  On the other hand, Chen and Deo's frequency domain test statistic is mean invariant
 }
\end{remark}

\begin{remark}
{\rm It seems natural to ask if a central limit theorem for statistics based on Bartlett's $T_p$ process can be obtained under the
GMC conditions on the errors. Although it might be possible to obtain a non-pivotal asymptotic null distribution under GMC conditions, the martingale transformation method  used in Delgado et al. (2005) and the frequency domain bootstrap approach in Hidalgo and Kreiss (2006)  may no longer be able to take care of the estimation effect  for the long memory model with unknown conditional heteroscedastic errors. The main reason is that the validity of both approaches rely on the assumption that the fourth order spectrum of the innovation sequence is a constant, which happens to be true for conditional homoscedastic martingale differences [cf. Shao (2010)]. In the case of conditional heteroscedastic errors, I am not aware of any feasible tests based on Bartlett's $T_p$ process. Further study along this direction would be certainly interesting.

% the asymptotic distribution of the $T_p$ process involves the fourth cumulant in the asymptotic covariance of the limiting process, and the validity of the martingale transformation method and the frequency domain approach  heavily depends on the assumption that the innovations are , which implies the constant fourth order cumulant spectrum [cf. Shao (2008)].

%Due the the involvement of the fourth cumulant effect,  be extend the  such results may not be of much practical interests. The main reason is that a detailed check of  the proofs of Delgado et al. (2005) and Hidalgo and Kreiss (2006) shows that  In the presence of conditional heteroscedastic errors, their methodologies  may no longer work.
%, and the vali resultseems to suggest that the fourth cumulant effect is not asymptotically negligible and  may not be able to remove the estimation effect under the assumptions that the innovations are conditional heteroscedastic.

}
\end{remark}

\section{Conclusions}
\label{sec:con}

In this paper, we showed that Hong's (1996) test is robust to conditional heteroscedasticity of unknown form in large sample theory
 and is  applicable to a large class of dependent white noise series. Further, when applied to the residuals from short/long memory time series models,
 the asymptotical null distribution is still valid.
   %It would be interesting to study the asymptotic behavior of the test statistic under local alternatives and  the techniques
 %presented here are expected to be useful. We leave this challenging task for future research.
 The main focus of this paper is on the theoretical aspect, although the empirical performance is also very important.  The finite sample performance of Hong's  test statistic has been examined by Hong (1996) and Chen and Deo (2004b) among others  to assess the goodness of fit of time series models with  iid errors. It was found that the sampling distribution of the test statistic is right-skewed, and the size distortion can presumably be reduced by adopting a power transformation method [Chen and Deo (2004b)] or frequency domain bootstrap approach [Paparoditis (2000)]. The performance of the afore-mentioned test statistics  along with size-correction devices have yet to be examined for time series models with dependent errors. An  in-depth study is certainly worthwhile, and will be pursued in a separate work.

\bigskip

 \baselineskip=17pt \centerline{REFERENCES}

\bigskip

%\par\noindent\hangindent2.3em\hangafter 1
%{ Anderson, T. W.} (1993). Goodness of fit tests for spectral
%distribution. {\it Annals of Statistics} { 21}, 830-847.

%\par\noindent\hangindent2.3em\hangafter 1
%{Bai, J.S.} (1993) On the partial sum processes of the residuals
%in autoregressive and moving average models. {\it Journal of Time
%Series Analysis}, 14, 247-260.

%\par\noindent\hangindent2.3em\hangafter 1
%{ Baillie, R. T.}, {C. F. Chung} \& {M. A. Tieslau} (1996).
%Analyzing inflation by the fractionally integrated ARFIMA-GARCH
%model.  {\it Journal of Applied Econometrics} { 11}, 23-40.

\par\noindent\hangindent2.3em\hangafter 1
{ Bartlett, M. S.} (1955). {\it An Introduction to Stochastic Processes
with Special Reference to Methods and Applications}. Cambridge
University Press.

\par\noindent\hangindent2.3em\hangafter 1
{ Beran, J.} (1992). A goodness-of-fit test for time series with
long range dependence. {\it Journal of  Royal Statistical Society  Series B Statistical Methodology} { 54}, 749-760.

\par\noindent\hangindent2.3em\hangafter 1
{ Beran, J.} (1995). Maximum likelihood  estimation of the
 differencing parameter for invertible short and long memory autoregressive integrated moving average models.
 {\it  Journal of  Royal Statistical Society  Series B Statistical Methodology} { 57}, 659-672.

%\par\noindent\hangindent2.3em\hangafter 1
%{ Bollerslev, T.} (1986). Generalized autoregressive conditional
%heteroscedasticity. {\it Journal of Econometrics} { 31}, 307-327.

\par\noindent\hangindent2.3em\hangafter 1
{ Box, G.} \& {D. Pierce} (1970). Distribution of residual
autocorrelations in autoregressive-integrated moving average time
series models. {\it Journal of the  American Statistical
Assocication} { 65}, 1509-1526.

\par\noindent\hangindent2.3em\hangafter 1
{ Breidt, F.J.}, {R. A. Davis} \& {A. A. Trindade} (2001).
 Least absolute deviation estimation for all-pass time series models. {\it
 Annals of   Statistics} { 29}, 919-946.

\par\noindent\hangindent2.3em\hangafter 1
{ Brillinger, D. R.} (1975). {\it Time Series: Data Analysis and
Theory}. Holden-Day, San Francisco.

\par\noindent\hangindent2.3em\hangafter 1
{ Chen, W.} \& {R. S. Deo} (2004a). A generalized portmanteau
goodness-of-fit test for time series models. {\it Econometric
Theory} { 20}, 382-416.

\par\noindent\hangindent2.3em\hangafter 1
{ Chen, W.} \& {R. S. Deo} (2004b). Power transformation to induce
normality and their applications. {\it Journal of  Royal Statistical Society  Series B Statistical Methodology} { 66}, 117-130.

\par\noindent\hangindent2.3em\hangafter 1
{ Chen, W.} \& {R. S. Deo} (2006). The variance ratio statistic at large horizons.
{\it Econometric Theory} 22, 206-234.

\par\noindent\hangindent2.3em\hangafter 1
{ Chen, H.} \& {J. P. Romano} (1999). Bootstrap-assisted
goodness-of-fit tests in the frequency domain. {\it Journal of
Time Series  Analysis} { 20}, 619-654.

%\par\noindent\hangindent2.3em\hangafter 1
%{ Dahlhaus, R.} (1989). Efficient parameter estimation for
%self-similar processes. {\it Annals of Statistics} { 17}, 1749-1766.

\par\noindent\hangindent2.3em\hangafter 1
{ Delgado, M.A.}, {J. Hidalgo} \& {C. Velasco} (2005).
Distribution free goodness-of-fit tests for linear processes. {\it
Annals of   Statistics} { 33}, 2568-2609.

%\par\noindent\hangindent2.3em\hangafter 1
%{ Delgado, M.A.} \& {C. Velasco} (2007). A new class of distribution-free tests
%for time series models specification. Preprint.

\par\noindent\hangindent2.3em\hangafter 1
{ Deo, R. S.} (2000). Spectral tests of the martingale hypothesis under
conditional heteroscedasticity. {\it Journal of Econometrics} { 99},
291-315.

%\par\noindent\hangindent2.3em\hangafter 1
%{ Ding, Z.}, {C. Granger} \& {R. Engle} (1993). A long memory property of stock market returns and
%a new model. {\it Journal of  Empirical Finance} { 1}, 83-106.

%\par\noindent\hangindent2.3em\hangafter 1
%{ Drouiche, K.} (2000). A new test for whiteness. {\it IEEE
%Transactions on Signal Processing} { 48}, 1864-1871.

%\par\noindent\hangindent2.3em\hangafter 1
%{ Drouiche, K.} (2007). A new test for spectrum flatness. {\it Journal of Time Series
%Analysis} { 28}, 793-806.

%\par\noindent\hangindent2.3em\hangafter 1
%{ Duchesne, P.} (2006). On testing for serial correlation with a
%wavelet-based spectral density estimator in multivariate time
%series. {\it Econometric Theory} { 22} 633-676.

\par\noindent\hangindent2.3em\hangafter 1
{ Durlauf, S.} (1991). Spectral based testing for the martingale
hypothesis. {\it Journal of Econometrics} { 50}, 1-19.

\par\noindent\hangindent2.3em\hangafter 1
{ Elek, P.} \& {L. M\'{a}rkus} (2004). A long range dependent
model with nonlinear innovations for simulating daily river flows.
{\it Natural Hazards and Earth System Sciences} 4, 277-283.

%\par\noindent\hangindent2.3em\hangafter 1
%{ Fan, J.}  \& {W. Zhang} (2004). Generalized likelihood ratio
%tests for spectral density. {\it Biometrika} { 91}, 195-209.

%\par\noindent\hangindent2.3em\hangafter 1
%{ Fox, R.} \& {M. S. Taqqu} (1986). Large-sample properties of
%parameter estimates for strongly dependent stationary Gaussian
%time series. {\it Annals of Statistics} { 14}, 517-532.

%\par\noindent\hangindent2.3em\hangafter 1
%{ Francq, C.} \& {J. M. Zako\"1an} (1998). Estimating
%linear representations of nonlinear processes. {\it Journal of
%Statistical  Planning and   Inference} { 68}, 145-165.

\par\noindent\hangindent2.3em\hangafter 1
{ Francq, C.} \& {J. M. Zako\"1an} (2000). Covariance
matrix estimation for estimators of mixing weak ARMA-models. {\it Journal of
 Statistical  Planning and  Inference} { 83}, 369-394.

\par\noindent\hangindent2.3em\hangafter 1
{ Francq, C.} \& {J. M. Zako\"1an} (2005). Recent Results
for Linear Time Series Models with Non Independent Innovations, in
Statistical Modeling and Analysis for Complex Data Problems, P.
Duchesne and B. Rémillard Editors,  Springer.

%\par\noindent\hangindent2.3em\hangafter 1
%{ Francq, C.} \& {J. M. Zako\"1an} (2006). HAC estimation
%and strong linearity testing in weak ARMA models, {\it Journal of
%Multivariate  Analysis} { 98}, 114-144.

\par\noindent\hangindent2.3em\hangafter 1
{ Francq, C.}, {R. Roy} \& {J. M. Zako\"1an} (2005).
Diagnostic Checking in ARMA models with uncorrelated errors. {\it
Journal of the  American  Statistical  Association.} { 100}, 532-544.

%\par\noindent\hangindent2.3em\hangafter 1
%{Giraitis, L.}, {Hidalgo, J.} and {P. M. Robinson} (2001) Gaussian
%estimation of parametric spectral density with unknown pole.  {\it
%Annals of Statistics} 29, 987-1023.

%\par\noindent\hangindent2.3em\hangafter 1
%{ Giraitis, L.} \& {D. Surgailis} (1990). A central limit theorem
%for quadratic forms in strongly dependent random variables and its
%application to asymptotic normality of Whittle's estimate. {\it
%Probability Theory and Related Fields} { 86}, 87-104.

\par\noindent\hangindent2.3em\hangafter 1
{ Granger, C. W. J.} \& {A. P. Anderson} (1978). {\it An Introduction to Bilinear Time Series Models.} Gottinger: Vandenhoek and Ruprecht.

\par\noindent\hangindent2.3em\hangafter 1
{ Granger, C. W. J.} \& {T. Ter\"asvirta} (1993). {\it Modelling Nonlinear Economic Relationships} (New York: Oxford University Press).

\par\noindent\hangindent2.3em\hangafter 1
{ Grenander, U.} \& {M. Rosenblatt} (1957). {\it Statistical
Analysis of Stationary Time Series}. Wiley, New York.

\par\noindent\hangindent2.3em\hangafter 1
{ Guerre, E.} \& {P. Lavergne} (2005) Data-driven rate-optimal specification testing in regression models.
{\it Annals of Statistics} 33, 840-870.

\par\noindent\hangindent2.3em\hangafter 1
{ Hall, P.} \& {C. C. Heyde} (1980). {\it Martingale Limit Theory and Its Applications}. Academic Press.

%\par\noindent\hangindent2.3em\hangafter 1
%{ Hannan, E. J.} (1973). The asymptotic theory of linear time series
%models. {\it Journal of  Applied Probability} { 10}, 130-145.

%\par\noindent\hangindent2.3em\hangafter 1
%{ Hauser, M. A.} \& {R. M. Kunst} (1998a). Fractionally integrated
%models with ARCH errors: with an application to the Swiss 1-month
%euromarket interest rate. {\it Review of Quantitative Finance and
%Accounting} { 10}, 95-113.

%\par\noindent\hangindent2.3em\hangafter 1
%{ Hauser, M. A.} \&  {R. M. Kunst} (1998b).  Forecasting
%high-frequency financial data with the ARFIMA-ARCH model. {\it
%Journal of Forecasting} { 20}, 501-518.

\par\noindent\hangindent2.3em\hangafter 1
{ Hidalgo, J.} \& {J. P. Kreiss} (2006). Bootstrap specification tests
for linear covariance stationary processes. {\it Journal of Econometrics} { 133}, 807-839.

\par\noindent\hangindent2.3em\hangafter 1
{ Hong, Y.} (1996). Consistent testing for serial correlation of
unknown form. {\it Econometrica} { 64}, 837-864.

\par\noindent\hangindent2.3em\hangafter 1
{ Hong, Y.} \& {Y. J. Lee} (2003). Consistent testing for serial
uncorrelation of unknown form under general conditional
heteroscedasticity. Preprint.

\par\noindent\hangindent2.3em\hangafter 1
{ Horowitz, J. L.}, {I. N. Lobato}, {J. C. Nankervis} \& {N. E. Savin}
(2006). Bootstrapping the Box-Pierce $Q$ test: A robust test of
uncorrelatedness. {\it Journal of Econometrics} { 133}, 841-862.

\par\noindent\hangindent2.3em\hangafter 1
{ Horowitz, J. L.} \& {V. G. Spokoiny} (2001) An adaptive rate-optimal test
of a parametric mean-regression model against a nonparametric alternative.
{\it Econometrica} 69, 599-631.

\par\noindent\hangindent2.3em\hangafter 1
{ Hosoya, Y.} (1997). A limit theory for long-range dependence and
statistical inference on related models.  {\it Annals of Statistics} { 25}, 105-137.

\par\noindent\hangindent2.3em\hangafter 1
{ Hsing, T.} \& {W. B. Wu} (2004). On weighted $U$-statistics for
stationary processes. {\it Annals of Probability}  { 32}, 1600-1631.

%\par\noindent\hangindent2.3em\hangafter 1
%{ Khmaladzee, V.} (1981). Martingale approach in the theory of goodness-of-fit tests. {\it Theory of
%Probability and its  Application} 26, 240-257.

\par\noindent\hangindent2.3em\hangafter 1
{ Koopman, S. J.}, {M. Oohs} \& {M. A. Carnero} (2007). Periodic
seasonal Reg-ARFIMA-GARCH models for daily electricity spot
prices. {\it Journal of the  American Statistical Association} { 102},
16-27.

%\par\noindent\hangindent2.3em\hangafter 1
%{ Lee, J.} \& {Y. Hong} (2001). Testing for serial correlation of
%unknown form using wavelet methods. {\it Econometric Theory} { 17},
%386-423.

\par\noindent\hangindent2.3em\hangafter 1
{ Li, G.} \& {W. K. Li} (2008) Least absolute deviation estimation for fractionally integrated autoregressive moving average time series models with conditional heteroscedasticity. {\it Biometrika} {95}, 399-414.

\par\noindent\hangindent2.3em\hangafter 1
{ Ling, S.} \& {W. K. Li} (1997) On fractionally integrated
autoregressive moving-average time series models with conditional
heteroscedasticity. {\it Journal of the American Statistical
Association} 92, 1184-1194.

%\par\noindent\hangindent2.3em\hangafter 1
%{Ljung, G. M.} and {Box, G. E. P.} (1978) On the measure of lack
%of fit in time series models. {\it Biometrika}, 65, 297-303.

\par\noindent\hangindent2.3em\hangafter 1
{ Lien, D.} \& {Y. K. Tse}  (1999). Forecasting the Nikkei spot
index with fractional cointegration. {\it Journal of Forecasting}
 { 18}, 259-273.

%\par\noindent\hangindent2.3em\hangafter 1
%{ Lobato, I.N.} (2001). Testing that a dependent process is
%uncorrelated. {\it Journal of the  American  Statistical
% Association} { 96}, 1066-1076.

\par\noindent\hangindent2.3em\hangafter 1
{ Lobato, I.N.}, {J. C. Nankervis} \& {N. E. Savin}  (2002). Testing
for zero autocorrelation in the presence of statistical
dependence. {\it Econometric Theory} { 18}, 730-743.

%\par\noindent\hangindent2.3em\hangafter 1
%{ Kl\"uppelberg, C.} and { Mikosch, T.} (1996). The integrated
%periodogram for stable processes. { 24}, 1855-1879.

\par\noindent\hangindent2.3em\hangafter 1
{ Lumsdaine, R. L.} \& {S. Ng} (1999). Testing for ARCH in the presence
of a possibly misspecified conditional mean.  {\it Journal of Econometrics} { 93}, 257-279.

%\par\noindent\hangindent2.3em\hangafter 1
%{ Mayoral, L.} (2007).  Minimum distance estimation of stationary
%and non-stationary ARFIMA processes.  {\it Econometrics Journal }
% { 10}, 124-148.

%\par\noindent\hangindent2.3em\hangafter 1
%{ Nordman, D. J.} \& {S. N. Lahiri} (2006). A frequency domain empirical likelihood
% for short- and long-range dependence. {\it Annals of Statistics} { 6}, 3019-3050.

\par\noindent\hangindent2.3em\hangafter 1
{ Paparoditis, E.} (2000). Spectral density based goodness-of-fit
tests for time series models, {\it Scandinavian   Journal of  Statistics} { 27}, 143-176.

\par\noindent\hangindent2.3em\hangafter 1
{ Priestley, M. B.} (1981). {\it Spectral Analysis and Time Series}, Vol 1,
Academic, New York.

\par\noindent\hangindent2.3em\hangafter 1
{ Robinson, P. M.} (2005). Efficiency improvements in inference on
stationary and nonstationary fractional time series. {\it Annals of
 Statistics} { 33}, 1800-1842.

\par\noindent\hangindent2.3em\hangafter 1
{ Romano, J. L.} \& {L. A. Thombs} (1996). Inference for
autocorrelations under weak assumptions. {\it Journal of the  American  Statistical  Association} {  91}, 590-600.

\par\noindent\hangindent2.3em\hangafter 1
{ Shao, X.}  (2010). Nonstationarity-extended Whittle estimation. {\it Econometric Theory}, to appear.

\par\noindent\hangindent2.3em\hangafter 1
{ Shao, X.} \& {W. B. Wu} (2007). Asymptotic spectral theory for
nonlinear time series. {\it Annals of Statistics} { 4}, 1773-1801.

%\par\noindent\hangindent2.3em\hangafter 1
%{ Subba Rao, T.} \& {M. M. Gabr} (1984). {\it An Introduction to
%Bispectral Analysis and Bilinear Time Series Models}. Lecture
%Notes in Statistics, 24. New York: Springer-Verlag.

%\par\noindent\hangindent2.3em\hangafter 1
%{ Tong, H.} (1990). {\it Non-linear Time Series: A Dynamical System
%Approach.} Oxford University Press.

%\par\noindent\hangindent2.3em\hangafter 1
%{ Velasco, C.} \& {P. M. Robinson} (2000). Whittle pseudo-maximum
%likelihood estimation for nonstationary time series. {\it Journal of the  American Statistical  Association} { 95}, 1229-1243.

\par\noindent\hangindent2.3em\hangafter 1
{ Wu, W. B.} (2005). Nonlinear system theory: another look at
dependence. {\it Proceedings of the National Academy of Science} { 102}, 14150-14154.

\par\noindent\hangindent2.3em\hangafter 1
{ Wu, W. B.} (2007).  Strong invariance principles for dependent
random variables. {\it Annals of Probability} { 35}, 2294-2320.

\par\noindent\hangindent2.3em\hangafter 1
{ Wu, W. B.} \& {W. Min} (2005).  On linear processes with
dependent innovations.  {\it Stochastic Processes and Their
Applications } { 115}, 939-958.

\par\noindent\hangindent2.3em\hangafter 1
{ Wu, W. B.} \& {X. Shao} (2004). Limit theorems for iterated
random functions. {\it Journal of   Applied Probability} { 41}, 425-436.

\par\noindent\hangindent2.3em\hangafter 1
{ Wu, W. B.} \& { X. Shao} (2007). A limit theorem for quadratic forms and its applications. {\it Econometric Theory} 23, 930-951.

\par\noindent\hangindent2.3em\hangafter 1
{ Wu, W. B.} \& {M. Woodroofe} (2004) Martingale approximations for sums of stationary processes. {\it
Annals of Probability} 32, 1674-1690.

\bigskip

\section{Technical Appendices}
\label{sec:tech}

Throughout the appendices, $u_t$ is assumed to be an uncorrelated stationary sequence with the representation (\ref{eq:ut}).
 For the convenience of notation,  let
$k_{nj}=K(j/m_n)$. Denote by
$Z_{jt}=u_tu_{t-j}$ and $D_{j,k}=\sum_{t=k}^{\infty}{\cal
P}_{k}(Z_{jt})$. Note that for each $j\in\N$, $D_{j,k}$ is a sequence of stationary and ergodic
martingale differences. For $a,b\in\R$, denote by
$a\vee b=\max(a,b)$ and $a\wedge b=\min(a,b)$. Let ${\cal
F}_i^j=(\varepsilon_i,\cdots,\varepsilon_j)$ and ${\cal
F}_t'=(\cdots,\varepsilon_{-1}',\varepsilon_0',\varepsilon_1,\cdots,\varepsilon_t)$,
$t\in\N$. For $X\in {\cal L}^1$, denote by ${\cal
P}_t'X=\E(X|{\cal F}_t')-\E(X|{\cal F}_{t-1}')$.
Let
$u_k^*=F(\cdots,\varepsilon_{-1},\varepsilon_0',\varepsilon_1,\cdots,\varepsilon_k)$, $k\in\N$.
Denote by $\delta_{\alpha}(k)=\|u_k-u_k^*\|_{\alpha}$, $k\in\N$, $\alpha\ge
1$ the physical dependence measure introduced by Wu (2005). According to  Wu
(2007), we have $\|{\cal P}_0Z_{jk}\|_{\alpha}\le
C(\delta_{2\alpha}(k)+\delta_{2\alpha}(k-j){\bf 1}(k\ge j))$ if $u_t\in {\cal L}^{2\alpha}$, and
$\delta_{\alpha}(k)\le Cr^k$ for some $r\in (0,1)$ provided that $u_t$ is  GMC$(\alpha)$, $\alpha\ge 1$.

One of major technical contributions of this paper is to replace the martingale difference assumption in Hong and Lee (2003) by the
GMC condition under the white noise null hypothesis. This is achieved by approximating the double array sequence $\sum_{t=j+1}^{n}Z_{jt}$ using its martingale counterpart $\sum_{t=j+1}^{n}D_{j,t}$ for  $j=1,\cdots,m_n$.
Note that the martingale approximation for the single array sequence $u_t$ has been well studied [cf. Hsing and Wu (2004), Wu and Woodroofe (2004), Wu and Shao (2007) among others], but the techniques there are not directly applicable. The major difficulty is that in our setting the martingale approximation error has to be bounded uniformly in $j=1,\cdots,m_n$ and the application of martingale central limit theorem after martingale approximation
 requires very delicate analysis due to the presence of dependence.

We separate the proofs of Theorem~\ref{th:null} and Theorem~\ref{th:farimaresidual} along with necessary lemmas into Appendices A and B respectively.
\subsection{Appendix A}

 Let $\theta_{j,r,\alpha}=\|{\cal P}_0 Z_{jr} \|_{\alpha}$,
 $\alpha\ge 1$ and
$\Theta_{j,n,\alpha}=\sum_{r=n}^{\infty}\theta_{j,r,\alpha}$. The
following lemma is an extension of Theorem 1 (ii) in Wu (2007).
Since the proof basically repeats that in Wu (2007), we omit the
details.
\begin{lemma}
\label{lem:diff} Assume $u_t\in {\cal L}^{2\alpha}$, $\alpha\ge
2$. For $0< a_n < b_n\le n$, we have
\[\left\|\sum_{r=a_n}^{b_n}(Z_{jr}-D_{jr})\right\|_{\alpha}^2\le C\sum_{k=1}^{b_n-a_n+1}\Theta_{j,k,\alpha}^2.\]
\end{lemma}

 The part (a) of the lemma below states the variance and covariances of the approximating martingale difference $D_{j,k}$ and may be of its independent interest.

\begin{lemma}
\label{lem:djk} Assume that $u_t$ is GMC$(8)$.
 (a) For $j>0$, we have
 \begin{eqnarray*}
 \E(D_{j,k}^2)=\sigma^4+\cov(u_t^2,u_{t-j}^2)+\sum_{k\not=0,k\in\Z}\cum(u_0,u_k,u_{-j},u_{k-j}),
 \end{eqnarray*}
and
$\E({D}_{j,k}{D}_{j',k})=(1/2)\sum_{k\in\Z}\{\cum(u_0,u_{-j},u_k,u_{k-j'})+\cum(u_0,u_{-j'},u_k,u_{k-j})\}$
when  $j\not=j'>0$.  (b) Let $D_{j,k}'=\sum_{t=k}^{\infty}{\cal P}_k'(u_t' u_{t-j}')$. Then  $\|D_{j,k}-D_{j,k}'\|_{4}\le
C(\rho^{k-j}{\bf 1}(k\ge j)+|j-k|{\bf 1}(k< j))$. (c) Let
$\tilde{D}_{j,k}=\E(D_{j,k}|(\varepsilon_k,\cdots,\varepsilon_{k-l+1}))$,
$l\in\N$. Then $\|\tilde{D}_{j,k}-D_{j,k}\|_{4}\le
C(\rho^{l-j}{\bf 1}(l\ge j)+|j-l|{\bf 1}(l<j))$. Here the positive constant $C$ appeared in (b) and (c) is independent of $j$.
\end{lemma}
 \noindent Proof of Lemma~\ref{lem:djk}: (a) It follows that when
 $j=j'>0$,
\begin{eqnarray*}
\E(D_{j,k}^2)&=&\sum_{k=-\infty}^{\infty}\cov(Z_{jt},Z_{j(t+k)})=\var(Z_{jt})+\sum_{k\not=0,
k\in\Z}\cov(u_{t}u_{t-j},u_{t+k}u_{t+k-j})\\
&=&\sigma^4+\cov(u_t^2,u_{t-j}^2)+\sum_{k\not=0,
k\in\Z}\cum(u_0,u_k,u_{-j},u_{k-j})
\end{eqnarray*}
 and when $j\not= j'>0$,
\begin{eqnarray*}
\E({D}_{j,k}{D}_{j',k})&=&(1/4)\E\{(D_{j,k}+D_{j',k})^2-(D_{j,k}-D_{j',k})^2\}\\
&=&(1/4)\sum_{k\in\Z}\{\cov(u_tu_{t-j}+u_tu_{t-j'},u_{t+k}u_{t+k-j}+u_{t+k}u_{t+k-j'})\\
&&-\cov(u_tu_{t-j}-u_tu_{t-j'},u_{t+k}u_{t+k-j}-u_{t+k}u_{t+k-j'})\}\\
&=&(1/2)\sum_{k\in\Z}\{\cov(u_tu_{t-j},u_{t+k}u_{t+k-j'})+\cov(u_tu_{t-j'},u_{t+k}u_{t+k-j})\}\\
&=&(1/2)\sum_{k\in\Z}\{\cum(u_0,u_{-j},u_k,u_{k-j'})+\cum(u_0,u_{-j'},u_k,u_{k-j})\}.
\end{eqnarray*}
(b) In general,  for $V_t=J(\cdots,\varepsilon_{t-1},\varepsilon_t)$, we have $\E(V_t|{\cal F}_k')=\E(V_t'|{\cal F}_k)$ when $t\ge k$. So for
$\alpha\ge 1$,
\begin{eqnarray*}
\|\E(V_t|{\cal F}_k)-\E(V_t'|{\cal F}_k')\|_{\alpha}&\le& \|\E(V_t|{\cal F}_k)-\E(V_t'|{\cal F}_k)\|_{\alpha}+\|\E(V_t|{\cal F}_k')-\E(V_t'|{\cal F}_k')\|_{\alpha}\\
&\le& 2\|V_t-V_{t}'\|_{\alpha},
\end{eqnarray*}
which implies that
\begin{eqnarray}
\label{eq:ineq2}
\|{\cal P}_k V_t-{\cal P}_k'V_{t'}\|_{\alpha}\le 4\|V_t-V_t'\|_{\alpha}.
\end{eqnarray}

Note that $D_{j,k}=\sum_{t=k}^{\infty}{\cal P}_k (u_tu_{t-j})$ and
$D_{j,k}'=\sum_{t=k}^{\infty}{\cal P}_k'(u_t'u_{t-j}')$. Then when
$k\le t\le k+j-1$, ${\cal P}_k (u_tu_{t-j})=u_{t-j}{\cal P}_k u_t$
and ${\cal P}_k'(u_t'u_{t-j}')=u_{t-j}'{\cal P}_k'u_t'$. So by the Cauchy-Schwarz inequality and (\ref{eq:ineq2}),
\begin{eqnarray*}
&&\hspace{-0.8cm}\|D_{j,k}-D_{j,k}'\|_{4}\le \sum_{t=k}^{k+j-1}\|u_{t-j}{\cal
P}_k u_t-u_{t-j}'{\cal
P}_k'u_t'\|_{4}+\sum_{t=k+j}^{\infty}\|{\cal P}_k
(u_tu_{t-j})-{\cal P}_k'(u_t'u_{t-j}')\|_{4}\\
&&\le C\sum_{t=k}^{k+j-1}\{\|u_{t-j}-u_{t-j}'\|_{8}+\|{\cal P}_k
u_t-{\cal
P}_{k}'u_{t}'\|_{8}\}+C\sum_{t=k+j}^{\infty}\|u_tu_{t-j}-u_t'u_{t-j}'\|_{4}\\
&&\le C\sum_{t=k}^{k+j-1}\{\rho^{t-j}+{\bf
1}(t\le j)+\rho^{t}\}+C\sum_{t=k+j}^{\infty}\{\rho^t+\rho^{t-j}\}\\
&&\le C\{\rho^{k-j}{\bf 1}(k\ge j)+|j-k|{\bf 1}(k< j)\}.
\end{eqnarray*}
As to (c), applying the fact that $\E(D_{j,l}|\varepsilon_l,\cdots,\varepsilon_1)=\E(D_{j,l}'|{\cal F}_l)$, we get
\begin{eqnarray*}
\|\tilde{D}_{j,k}-D_{j,k}\|_{4}&=&\|\tilde{D}_{j,l}-D_{j,l}\|_{4}=\|D_{j,l}-\E(D_{j,l}|\varepsilon_l,\cdots,\varepsilon_1)\|_{4}\\
&=&\|\E((D_{j,l}-D_{j,l}')|{\cal F}_{l})\|_{4}\le
\|D_{j,l}-D_{j,l}'\|_{4}\\
&\le& C\{\rho^{l-j}{\bf 1}(l\ge j)+|j-l|{\bf 1}(l<j)\}.
\end{eqnarray*}
The proof is complete. \qed

\bigskip

 {\it Proof of Theorem~\ref{th:null}}:  Since
$\hat{R}_{u}(0)=\sigma^2+O_p(n^{-1/2})$, we have
\[n\sum_{j=1}^{m_n}k_{nj}^2\hat{\rho}_{u}^2(j)=n\sigma^{-4}\sum_{j=1}^{m_n}k_{nj}^2\hat{R}_u^2(j)+o_p(m_n^{1/2}).\]
 Let
$G_n:=n\sum_{j=1}^{m_n}k_{nj}^2\tilde{R}_{u}^2(j)$, where
$\tilde{R}_u(j)=n^{-1}\sum_{t=|j|+1}^{n}u_t u_{t-|j|}$.
Note that
$\tilde{R}_u(j)-\hat{R}_u(j)=\bar{u}\{(1-j/n)\bar{u}-n^{-1}\sum_{t=1}^{n-j+1}u_t-n^{-1}\sum_{t=j+1}^{n}u_t\}$ for $j\ge 1$.
Under GMC$(2)$,  $\bar{u}^2=O_p(n^{-1})$,  $\sum_{j=1}^{m_n}k_{nj}^2\E(\sum_{t=j+1}^{n}u_t+\sum_{t=1}^{n-j+1}u_t)^2=O(n m_n)$.
Consequently,
$n\sum_{j=1}^{m_n}k_{nj}^2(\tilde{R}_u(j)-\hat{R}_u(j))^2=o_p(1)$.
Then it suffices to show
\begin{eqnarray}
\label{eq:Gn} \frac{G_n-\sigma^4 m_n C(K)}{(2\sigma^8 m_n
D(K))^{1/2}}\rightarrow_{D} N(0,1).
\end{eqnarray}
We shall approximate $G_n$ by
$\tilde{G}_n=\sum_{j=1}^{m_n}k_{nj}^2n^{-1}\left(\sum_{k=j+1}^{n}D_{j,k}\right)^2$.
  By the Cauchy-Schwarz inequality,
 \begin{eqnarray*}
|G_n-\tilde{G}_n|^2
%&=&\sum_{j=1}^{m_n}K_{nj}^2n^{-2}\left(\sum_{k=j+1}^{n}(Z_{jk}-D_{j,k})\right)\left(\sum_{k=j+1}^{n}(Z_{jk}+D_{j,k})\right)\\
          &\le&\sum_{j=1}^{m_n}\frac{k_{nj}^2}{n}\left(\sum_{k=j+1}^{n}(Z_{jk}-D_{j,k})\right)^2\times
          \sum_{j=1}^{m_n}\frac{k_{nj}^2}{n}\left(\sum_{k=j+1}^{n}(Z_{jk}+D_{j,k})\right)^2,
      %    &=&O_p(m_n/n)O_p(m_n)=o_p(m_n).
 \end{eqnarray*}
 where the second term on the right hand side of the inequality is easily shown to be $O_p(m_n)$ in view of the proof to be presented hereafter. As to the first term, we apply Lemma~\ref{lem:diff} and get
 \begin{eqnarray*}
 \frac{1}{n}\sum_{j=1}^{m_n}\left\|\sum_{k=j+1}^{n}(Z_{jk}-D_{j,k})\right\|^2&\le&\frac{C}{n}\sum_{j=1}^{m_n}\sum_{h=1}^{\infty}\left(\sum_{k=h}^{\infty}\|{\cal P}_0
Z_{jk}\|\right)^2\\
&\le&\frac{C}{n}\sum_{j=1}^{m_n}\sum_{h=1}^{\infty}\left(\sum_{k=h}^{\infty}(\delta_4(k)+\delta_4(k-j){\bf
1}(k\ge j))\right)^2\\
&\le&\frac{C}{n}\sum_{j=1}^{m_n}\sum_{h=1}^{\infty}\left(\sum_{k=h}^{\infty}(\delta_4(k)+\delta_4(k-j){\bf
1}(k\ge j))\right)\\
&\le&\frac{Cm_n}{n}\sum_{k=1}^{\infty}k\delta_4(k)+\frac{C}{n}\sum_{k=1}^{\infty}\sum_{h=1}^{k}\sum_{j=1}^{m_n\wedge
k}\delta_4(k-j)\\
&\le&Cm_n^2/n=o(1).
\end{eqnarray*}
So $G_n=\tilde{G}_n+o_p(m_n^{1/2})$.
Write
\begin{eqnarray*}
\tilde{G}_n&=&n^{-1}\sum_{j=1}^{m_n}k_{nj}^2\left(\sum_{k=j+1}^{n}D_{j,k}\right)^2\\
         &=&n^{-1}\sum_{j=1}^{m_n}k_{nj}^2\sum_{k=j+1}^{n}D_{j,k}^2+2n^{-1}\sum_{j=1}^{m_n}k_{nj}^2\sum_{k=j+2}^n\sum_{r=j+1}^{k-1}D_{j,k}D_{j,r}=\tilde{G}_{1n}+\tilde{G}_{2n}.
\end{eqnarray*}
Under the assumption that $u_t$ is GMC$(8)$, it is easy to show that $u_t^2$ is GMC$(4)$, which
implies that $|\cov(u_t^2,u_{t-j}^2)|\le Cr^j$ for some $r\in
(0,1)$. So by Lemma~\ref{lem:djk},
\begin{eqnarray*}
\E(\tilde{G}_{1n})
%&=&n^{-1}\sum_{j=1}^{m_n}k_{nj}^2(n-j)\E(D_{j,k}^2)\\
                 &=&n^{-1}\sum_{j=1}^{m_n}k_{nj}^2(n-j)\left(\sigma^4+\cov(u_t^2,u_{t-j}^2)+\sum_{k\not=0}\cum(u_0,u_k,u_{-j},u_{k-j})\right)\\
                 &=&\sigma^4\sum_{j=1}^{m_n}k_{nj}^2+O(1)=\sigma^4 m_n
                 C(K)+O(1),
\end{eqnarray*}
where we have applied  the absolute summability of the   $4$-th joint cumulants
under GMC$(4)$ [Wu and Shao (2004), Proposition 2].
 Let
$\tilde{D}_{j,k}=\E(D_{j,k}|\varepsilon_k,\varepsilon_{k-1},\cdots,\varepsilon_{k-l+1})$,
where $l=l_n=2m_n$. By Lemma~\ref{lem:djk} and the assumption that
$\log n=o(m_n)$,
\begin{eqnarray}
\label{eq:bound1} \sup_{1\le j\le
m_n}\|\tilde{D}_{j,k}-D_{j,k}\|_4=O(n^{-\kappa})~\mbox{for
any}~\kappa>0.
\end{eqnarray}
Write
\begin{eqnarray*}
\tilde{G}_{1n}&=&n^{-1}\sum_{j=1}^{m_n}k_{nj}^2\sum_{k=j+1}^{n}\tilde{D}_{j,k}^2+n^{-1}\sum_{j=1}^{m_n}k_{nj}^2\sum_{k=j+1}^{n}(D_{j,k}^2-\tilde{D}_{j,k}^2)=\tilde{G}_{11n}+\tilde{G}_{12n},
\end{eqnarray*}
where $\var(\tilde{G}_{11n})=O(m_n^3/n)=o(m_n)$ by the
$l_n$-dependence of $\tilde{D}_{j,k}$, and by (\ref{eq:bound1}),
\[\|\tilde{G}_{12n}\|\le \frac{C}{n}\sum_{j=1}^{m_n}k_{nj}^2\sum_{k=j+1}^{n}\|D_{j,k}^2-\tilde{D}_{j,k}^2\|=o(1).\]
 So (\ref{eq:Gn}) follows if we can show that
$\tilde{G}_{2n}/(2\sigma^8 m_n D(K))^{1/2}\rightarrow_{D}N(0,1)$.

Write
\begin{eqnarray}
\label{eq:tildeG2n}
\tilde{G}_{2n}&=&2n^{-1}\sum_{j=1}^{m_n}k_{nj}^2\times
\left(\sum_{k=j+2}^{6m_n}\sum_{r=j+1}^{k-1}+\sum_{k=6m_n+1}^{n}\sum_{r=j+1}^{m_n+1}+\sum_{k=6m_n+1}^{n}\sum_{r=k-2l_n+1}^{k-1}\right.\nonumber\\
&&\left.+\sum_{k=6m_n+1}^{n}\sum_{r=m_n+2}^{k-2l_n}\right){D}_{j,k}{D}_{j,r}=U_{1n}+U_{2n}+U_{3n}+U_{4n}.
\end{eqnarray}
We proceed to show that $U_{kn}=o_p(m_n^{1/2})$, $k=1,2,3$. Note
that the summands in $U_{1n}$ form martingale differences. So
\begin{eqnarray*}
\E(U_{1n}^2)=
\frac{4}{n^2}\sum_{k=3}^{6m_n}\left\|\sum_{r=2}^{k-1}k_{nj}^2\sum_{j=1}^{(r-1)\wedge
m_n}{D}_{j,k}{D}_{j,r}\right\|^2=O(m_n^5/n^2)=o(m_n).
\end{eqnarray*}
Regarding $U_{2n}$, we let
$\tilde{U}_{2n}=2n^{-1}\sum_{j=1}^{m_n}k_{nj}^2\sum_{k=6m_n+1}^{n}\sum_{r=j+1}^{m_n+1}\tilde{D}_{j,k}\tilde{D}_{j,r}$.
It is easy to show that $U_{2n}-\tilde{U}_{2n}=o_p(1)$ in view of
(\ref{eq:bound1}). Further, by Lemma~\ref{lem:djk},

 \begin{eqnarray*}
\E(\tilde{U}_{2n}^2)
%&=&\frac{4}{n^2}\sum_{k,k'=6m_n+1}^{n}\sum_{j,j'=1}^{m_n}k_{nj}^2k_{nj'}^2\sum_{r=j+1}^{m_n+1}\sum_{r'=j'+1}^{m_n+1}\E(\tilde{D}_{j,k}\tilde{D}_{j,r}\tilde{D}_{j',k'}\tilde{D}_{j',r'})\\
            &=&\frac{4}{n^2}\sum_{k,k'=6m_n+1}^{n}\sum_{j,j'=1}^{m_n}k_{nj}^2k_{nj'}^2\sum_{r=j+1}^{m_n+1}\sum_{r'=j'+1}^{m_n+1}\E(\tilde{D}_{j,k}\tilde{D}_{j',k'})\E(\tilde{D}_{j,r}\tilde{D}_{j',r'})\\
            &=&\frac{4(1+o(1))}{n^2}\sum_{k=6m_n+1}^{n}\sum_{j,j'=1}^{m_n}k_{nj}^2k_{nj'}^2\sum_{r=(j+1)\vee(j'+1)}^{m_n+1}\E(\tilde{D}_{j,k}\tilde{D}_{j',k})\E(\tilde{D}_{j,r}\tilde{D}_{j',r})\\
            &=&O(m_n^3/n)=o(m_n).
 \end{eqnarray*}
Thus $U_{2n}=o_p(m_n^{1/2})$. Concerning $U_{3n}$, since it is a
martingale, we have
 \begin{eqnarray*}
\E(U_{3n}^2)&=&\frac{4}{n^2}\sum_{k=6m_n+1}^{n}\left\|\sum_{j=1}^{m_n}k_{nj}^2\sum_{r=k-2l_n+1}^{k-1}D_{j,k}D_{j,r}\right\|^2\\
&\le&\frac{C}{n^2}\sum_{k=6m_n+1}^{n}\left(\sum_{j=1}^{m_n}\left\|\sum_{r=k-2l_n+1}^{k-1}D_{j,k}D_{j,r}\right\|\right)^2\\
&\le&\frac{C}{n^2}\sum_{k=6m_n+1}^{n}\left(\sum_{j=1}^{m_n}\left\|\sum_{r=k-2l_n+1}^{k-1}D_{j,r}\right\|_4\right)^2.
\end{eqnarray*}
Since $D_{j,r}$'s are martingale differences for each $j$, we
apply Burkholder's inequality [Hall and Heyde (1980)] and get
\[\left\|\sum_{r=k-2l_n+1}^{k-1}D_{j,r}\right\|_4\le C\left\|\sum_{r=k-2l_n+1}^{k-1}D_{j,r}^2\right\|^{1/2}\le C \left(\sum_{r=k-2l_n+1}^{k-1}\|D_{j,r}^2\|\right)^{1/2} \le Cm_n^{1/2}.\]
Note that the constant $C$ in the above display does not depend on $j$.
So $\E(U_{3n}^2)\le Cm_n^3/n=o(m_n)$. Let
$\tilde{U}_{4n}=2n^{-1}\sum_{j=1}^{m_n}k_{nj}^2\sum_{k=6m_n+1}^{n}\sum_{r=m_n+2}^{k-2l_n}\tilde{D}_{j,k}\tilde{D}_{j,r}$.
Since $U_{4n}-\tilde{U}_{4n}=o_p(1)$ by (\ref{eq:bound1}),   it
remains to show $\tilde{U}_{4n}/(2\sigma^8
m_nD(K))^{1/2}\rightarrow_{D}N(0,1)$ in view of
(\ref{eq:tildeG2n}).

Write $\tilde{U}_{4n}=n^{-1}\sum_{k=6m_n+1}^{n}V_{nk}$, where
$V_{nk}:=2\sum_{r=m_n+2}^{k-2l_n}\sum_{j=1}^{m_n}k_{nj}^2\tilde{D}_{j,k}\tilde{D}_{j,r}$.
Then $\{V_{nk}\}$ forms a sequence of martingale differences with
respect to ${\cal F}_{k}$.
 By the martingale central limit
theorem, it suffices to verify the following conditions:
\begin{eqnarray}
\label{eq:var}
\sigma^2(n):=\E(\tilde{U}_{4n}^2)&=&2\sigma^8 m_n D(K)(1+o(1)),\\
\label{eq:linderberg} \sum_{t=6m_n+1}^{n}\E(V_{nt}^2{\bf
1}(|V_{nt}|>\epsilon n\sigma(n)))&=&o(\sigma^2(n)n^2),~\epsilon>0,\\
\label{eq:conditionalvar}
\sigma^{-2}(n)n^{-2}\sum_{t=6m_n+1}^n\bar{V}_{nt}^2&\rightarrow_{p}&
1,~\mbox{where}~\bar{V}_{nt}^2=\E(V_{nt}^2|{\cal F}_{t-1}).
\end{eqnarray}
By Lemma~\ref{lem:djk} and (\ref{eq:bound1}), we have
\begin{eqnarray}
\label{eq:sigman}
\sigma^2(n)&=&n^{-2}\sum_{k=6m_n+1}^{n}\E(V_{nk}^2)\nonumber\\
           &=&\frac{4}{n^{2}}\sum_{k=6m_n+1}^{n}\sum_{r,r'=m_n+2}^{k-2l_n}\sum_{j,j'=1}^{m_n}k_{nj}^2k_{nj'}^2\E(\tilde{D}_{j,k}
           \tilde{D}_{j',k})\E(\tilde{D}_{j,r}\tilde{D}_{j',r'})\\
           &=&\frac{4}{n^{2}}\sum_{k=6m_n+1}^{n}\sum_{r=m_n+2}^{k-2l_n}\sum_{j,j'=1}^{m_n}k_{nj}^2k_{nj'}^2\E(\tilde{D}_{j,k}
           \tilde{D}_{j',k})\E(\tilde{D}_{j,r}\tilde{D}_{j',r})\nonumber\\
           &=&\frac{4}{n^{2}}\sum_{k=6m_n+1}^{n}\sum_{r=m_n+2}^{k-2l_n}\sum_{j,j'=1}^{m_n}k_{nj}^2k_{nj'}^2\E({D}_{j,k}
           {D}_{j',k})\E({D}_{j,r}{D}_{j',r})+o(1)\nonumber\\
           &=&\frac{4}{n^{2}}\sum_{k=6m_n+1}^{n}\sum_{r=m_n+2}^{k-2l_n}\sum_{j=1}^{m_n}k_{nj}^4\E(D_{j,k}^2)\E(D_{j,r}^2)(1+o(1))\nonumber\\
           &=&2\sigma^8m_nD(K)+o(m_n).\nonumber
 \end{eqnarray}
For (\ref{eq:linderberg}), again by Burkholder's inequality, we
get
\begin{eqnarray*}
\label{eq:cltcond1}
\E(V_{nk}^4)&=&\E\left(\sum_{r=m_n+2}^{k-2l_n}\sum_{j=1}^{m_n}k_{nj}^2\tilde{D}_{j,k}\tilde{D}_{j,r}\right)^4\le Cm_n^3\sum_{j=1}^{m_n}\E\left(\sum_{r=m_n+2}^{k-2l_n}\tilde{D}_{j,k}\tilde{D}_{j,r}\right)^4\\
&\le&Cm_n^3\sum_{j=1}^{m_n}\E(\tilde{D}_{j,k}^4)\E\left(\sum_{r=m_n+2}^{k-2l_n}\tilde{D}_{j,r}^2\right)^2\le
Cm_n^4 k^2,
\end{eqnarray*}
which implies (\ref{eq:linderberg}).  To show
(\ref{eq:conditionalvar}), we let
$\bar{V}_n^2=n^{-2}\sum_{t=6m_n+1}^{n}\bar{V}_{nt}^2$, where
\begin{eqnarray*}
\bar{V}_{nt}^2
%=\E\left(\left(\sum_{r=m_n+2}^{t-2l_n}\sum_{j=1}^{m_n}k_{nj}^2\tilde{D}_{j,t}\tilde{D}_{j,r}\right)^2\left|\right.{\cal
%F}_{t-1}\right)\\
&=&4\sum_{r,r'=m_n+2}^{t-2l_n}\sum_{j,j'=1}^{m_n}k_{nj}^2k_{nj'}^2\E(\tilde{D}_{j,t}\tilde{D}_{j',t}|{\cal
F}_{t-1})\tilde{D}_{j,r}\tilde{D}_{j',r'}.
\end{eqnarray*}
 Then we can write
 \begin{eqnarray}
 \label{eq:vnbar}
 \bar{V}_n^2-\sigma^2(n)&=&\frac{4}{n^2}\sum_{t=6m_n+1}^{n}\sum_{r,r'=m_n+2}^{t-2l_n}\sum_{j,j'=1}^{m_n}k_{nj}^2k_{nj'}^2\\
 &&\{\E([\tilde{D}_{j,t}\tilde{D}_{j',t}-{D}_{j,t}{D}_{j',t}]|{\cal
F}_{t-1})\tilde{D}_{j,r}\tilde{D}_{j',r'}\nonumber\\
&&+[\E({D}_{j,t}{D}_{j',t}|{\cal
F}_{t-1})-\E({D}_{j,t}{D}_{j',t}|{\cal F}_{t-l+1}^{t-1})]\tilde{D}_{j,r}\tilde{D}_{j',r'}\nonumber\\
&&+[\E({D}_{j,t}{D}_{j',t}|{\cal
F}_{t-l+1}^{t-1})-\E({D}_{j,t}{D}_{j',t})]
\tilde{D}_{j,r}\tilde{D}_{j',r'}\nonumber\\
&&+\E({D}_{j,t}{D}_{j',t})[\tilde{D}_{j,r}\tilde{D}_{j',r'}-\E(\tilde{D}_{j,r}\tilde{D}_{j',r'})]\nonumber\\
&&+\E({D}_{j,t}{D}_{j',t})\E(\tilde{D}_{j,r}\tilde{D}_{j',r'})\}-\sigma^2(n)=:\sum_{k=1}^{5}J_{kn}-\sigma^2(n).\nonumber
\end{eqnarray}

By a similar argument as in (\ref{eq:sigman}),  $J_{5n}=\sigma^2(n)(1+o(1))$. So (\ref{eq:conditionalvar}) follows if we can
show $\sigma^{-2}(n)J_{kn}=o_p(1)$ for $k=1,\cdots,4$. By
(\ref{eq:bound1}),  $J_{1n}=o_p(m_n)$. As to $J_{2n}$, it follows
from Lemma~\ref{lem:djk} and (\ref{eq:bound1}) that uniformly in $j,j'=1,2,\cdots,m_n$,
\begin{eqnarray*}
&&\hspace{-0.5cm}\|\E({D}_{j,t}{D}_{j',t}|{\cal
F}_{t-1})-\E({D}_{j,t}{D}_{j',t}|{\cal
F}_{t-l+1}^{t-1})\|=\|\E({D}_{j,l}{D}_{j',l}|{\cal
F}_{l-1})-\E({D}_{j,l}{D}_{j',l}|{\cal F}_{1}^{l-1})\|\\
&&\hspace{0.5cm}\le\|\E(({D}_{j,l}{D}_{j',l}-{D}_{j,l}'{D}_{j',l}')|{\cal
F}_{l-1})\|+\|\E(({D}_{j,l}{D}_{j',l}-{D}_{j,l}'{D}_{j',l}')|{\cal
F}_1^{l-1})\|\\
&&\hspace{0.5cm}\le
2\|{D}_{j,l}{D}_{j',l}-{D}_{j,l}'{D}_{j',l}'\|\le C\rho^{m_n}=O(n^{-\kappa})~\mbox{for any}~\kappa>0.
\end{eqnarray*}
So $J_{2n}=o_p(m_n)$.
Lemmas~\ref{lem:J3n} and ~\ref{lem:J4n} assert that $J_{3n}=o_p(m_n)$ and $J_{4n}=o_p(m_n)$ respectively.
Thus (\ref{eq:conditionalvar}) holds and the conclusion follows.
\qed

\begin{lemma}
\label{lem:J3n} Under the assumptions in Theorem~\ref{th:null},
the random variable $J_{3n}=4/n^2\sum_{t=6m_n+1}^{n}\sum_{r,r'=m_n+2}^{t-2l_n}\sum_{j,j'=1}^{m_n}k_{nj}^2k_{nj'}^2[\E({D}_{j,t}{D}_{j',t}|{\cal
F}_{t-l+1}^{t-1})-\E({D}_{j,t}{D}_{j',t})]
\tilde{D}_{j,r}\tilde{D}_{j',r'}$ as defined in (\ref{eq:vnbar}) is $o_p(m_n)$.
\end{lemma}
\noindent Proof of Lemma~\ref{lem:J3n}: Let
$M(j,j';t)=\E({D}_{j,t}{D}_{j',t}|{\cal
F}_{t-l+1}^{t-1})-\E({D}_{j,t}{D}_{j',t})$ and
\[\tilde{J}_{3n}=\frac{4}{n^2}\sum_{t=6m_n+1}^{n}\sum_{r,r'=m_n+2}^{t-2l_n}\sum_{j,j'=1}^{m_n}k_{nj}^2k_{nj'}^2M(j,j';t){D}_{j,r}{D}_{j',r'}.\]
It is easy to see that $\tilde{J}_{3n}=J_{3n}+o_p(m_n)$ in view of
(\ref{eq:bound1}). For notational convenience,  denote by
$H_D(j,t)=\sum_{r=m_n+2}^{t-2l_n}D_{j,r}$ and
$H_Z(j,t)=\sum_{r=m_n+2}^{t-2l_n}Z_{jr}$.
Write $\tilde{J}_{3n}={J}_{31n}+{J}_{32n}+J_{33n}$, where
\begin{eqnarray*}
{J}_{31n}&=&\frac{4}{n^2}\sum_{t=6m_n+1}^{n}\sum_{j,j'=1}^{m_n}k_{nj}^2k_{nj'}^2M(j,j';t)(H_D(j,t)-H_Z(j,t))
H_D(j',t),\\
{J}_{32n}&=&\frac{4}{n^2}\sum_{t=6m_n+1}^{n}\sum_{j,j'=1}^{m_n}k_{nj}^2k_{nj'}^2M(j,j';t)H_Z(j,t)
(H_D(j',t)-H_Z(j',t)), \\
J_{33n}&=&\frac{4}{n^2}\sum_{t=6m_n+1}^{n}\sum_{j,j'=1}^{m_n}k_{nj}^2k_{nj'}^2M(j,j';t)H_Z(j,t)H_Z(j',t).
\end{eqnarray*}
 We shall first prove $J_{31n}=o_p(m_n)$.
 Since $M(j,j',t)$ is $l_n$-dependent with respect $t$, we get by
  the Cauchy-Schwarz inequality,
\begin{eqnarray*}
\E({J}_{31n}^2)&\le&\frac{C}{n^4}\sum_{t=6m_n+1}^{n}\sum_{t'=(6m_n+1)\vee
(t-l_n)}^{n\wedge
(t+l_n)}\sum_{j_1,j_1',j_2,j_2'=1}^{m_n}\left\|H_{D}(j_1,t)-H_{Z}(j_1,t)\right\|_4\\
&&\left\|H_{D}(j_2,t')-H_{Z}(j_2,t')\right\|_4\left\|H_D(j_1',t)\right\|_4
\left\|H_D(j_2',t')\right\|_4.
\end{eqnarray*}
Since the summands in $H_D(j,t)$ form martingale differences, we apply
Burkholder's inequality and obtain
\begin{eqnarray}
\label{eq:burk} \|H_D(j,t)\|_4^4\le
C\E\left(\sum_{r=m_n+1}^{t-2l_n}D_{j,r}^2\right)^2\le
Ct^{2},~j=1,2,\cdots,m_n.
\end{eqnarray}
 Applying Lemma~\ref{lem:diff} and the fact that $\delta_8(k)\le Cr^k$ for some
$r\in (0,1)$, we get
\begin{eqnarray}
\label{eq:burk2}
&&\sum_{j_1=1}^{m_n}\left\|H_{D}(j_1,t)-H_{Z}(j_1,t)\right\|_4\le
C\sum_{j_1=1}^{m_n}
\left(\sum_{k_1=1}^{t-5m_n-1}\Theta_{j_1,k_1,4}^2\right)^{1/2}\nonumber\\
%&&\le C\sum_{j_1=1}^{m_n}
%\left(\sum_{k_1=1}^{t-5m_n-1}\Theta_{j_1,k_1,4}\right)^{1/2}\nonumber\\
&&\le C\sum_{j_1=1}^{m_n}\left(\sum_{k_1=1}^{t-5m_n-1}\sum_{h=k_1}^{\infty}(\delta_8(h)+\delta_8(h-j_1){\bf
1}(h\ge j_1))\right)^{1/2}\le Cm_n^{3/2}.
\end{eqnarray}
 Therefore, in view of (\ref{eq:burk}) and (\ref{eq:burk2}), we
 obtain $\E({J}_{31n}^2)\le Cm_n^6/n^2=o(m_n^2)$.
To show $J_{32n}=o_p(m_n)$, we note that
\begin{eqnarray}
\label{eq:expansion}
&&\hspace{-0.5cm}\|H_{Z}(j,t)\|_4^4=\sum_{r_1,r_2,r_2,r_4=m_n+1}^{t-2l_n}\E(Z_{jr_1}Z_{jr_2}Z_{jr_3}Z_{jr_4})\nonumber\\
             &&=\sum_{r_1,r_2,r_2,r_4=m_n+1}^{t-2l_n}\{\cov(Z_{jr_1},Z_{jr_2})\cov(Z_{jr_3},Z_{jr_4})+\cov(Z_{jr_1},Z_{jr_3})\cov(Z_{jr_2},Z_{jr_4})\nonumber\\
             &&+\cov(Z_{jr_1},Z_{jr_4})\cov(Z_{jr_2},Z_{jr_3})+\cum(Z_{jr_1},Z_{jr_2},Z_{jr_3},Z_{jr_4})\}.
\end{eqnarray}
Since $\{u_t\}$ are uncorrelated and the $k$-th ($k=2,3,\cdots,8$)
joint cumulants are absolutely summable under GMC$(8)$ [see Wu and
Shao (2004) Proposition 2], it is not hard to see that
$\|H_{Z}(j,t)\|_4^4\le Ct^2$. Following the same argument as in
the derivation of $\E(J_{31n}^2)$, we can derive
$\E(J_{32n}^2)=o(m_n^2)$, so $J_{32n}=o_p(m_n)$.

It remains to show that  ${J}_{33n}=o_p(m_n)$. Note that
\begin{eqnarray*}
\E({J}_{33n}^2)&\le&\frac{C}{n^4}\sum_{t=6m_n+1}^{n}\sum_{t'=(6m_n+1)\vee
(t'-l_n)}^{n\wedge
(t+l_n)}\sum_{j_1,j_1',j_2,j_2'=1}^{m_n}\sum_{r_1,r_2=m_n+2}^{t-2l_n}\sum_{r_1',r_2'=m_n+2}^{t'-2l_n}\\
&&|\E(Z_{j_1r_1}Z_{j_2r_2}Z_{j_1'r_1'}Z_{j_2'r_2'})|\le
\frac{C}{n^4}\sum_{t=6m_n+1}^{n}\sum_{t'=(6m_n+1)\vee
(t'-l_n)}^{n\wedge (t+l_n)}H_n(t,t').
\end{eqnarray*}
Following (\ref{eq:expansion}),  we can write $\E(Z_{j_1r_1}Z_{j_2r_2}Z_{j_1'r_1'}Z_{j_2'r_2'})$
as a sum of four components, which implies
%\begin{eqnarray*}
%&&\hspace{-0.5cm}\E(Z_{j_1r_1}Z_{j_2r_2}Z_{j_1'r_1'}Z_{j_2'r_2'})=\cov(Z_{j_1r_1},Z_{j_2r_2})\cov(Z_{j_1'r_1'},Z_{j_2'r_2'})+\cov(Z_{j_1r_1},Z_{j_1'r_1'})\\
%&&\hspace{-0.5cm}\cov(Z_{j_2r_2},Z_{j_2'r_2'})+\cov(Z_{j_1r_1},Z_{j_2'r_2'})\cov(Z_{j_2r_2},Z_{j_1'r_1'})+\cum(Z_{j_1r_1},Z_{j_2r_2},Z_{j_1'r_1'},Z_{j_2'r_2'}),
%\end{eqnarray*}
 $H_n(t,t')=\sum_{k=1}^{4}H_{kn}(t,t')$. For
$H_{1n}(t,t')$, it follows from the absolute summability of the
$4$-th cumulant that
\begin{eqnarray*}
H_{1n}(t,t')&=&\sum_{j_1,j_1',j_2,j_2'=1}^{m_n}\sum_{r_1,r_2=m_n+2}^{t-2l_n}\sum_{r_1',r_2'=m_n+2}^{t'-2l_n}|\{\cov(u_{r_1},u_{r_2})\cov(u_{r_1-j_1},u_{r_2-j_2})\\
&&+\cum(u_{r_1},u_{r_1-j_1},u_{r_2},u_{r_2-j_2})\}\{\cov(u_{r_1'},u_{r_2'})\cov(u_{r_1'-j_1'},u_{r_2'-j_2'})\\
&&+\cum(u_{r_1'},u_{r_2'},u_{r_1'-j_1'},u_{r_2'-j_2'})\}|\le
Cm_n^2(t \vee t')^2.
\end{eqnarray*}
By the same argument, we have $H_{kn}(t,t')\le Cm_n^2 (t\vee
t')^2$, $k=2,3$. Regarding $H_{4n}(t,t')$, we apply the product
theorem for the joint cumulants [Brillinger (1975)] and write
\begin{eqnarray*}
\cum(Z_{j_1r_1},Z_{j_2r_2},Z_{j_1'r_1'},Z_{j_2'r_2'})=\sum_{v}\cum(u_{i_j},i_j\in
v_1)\cdots\cum(u_{i_j},i_j\in v_p),
\end{eqnarray*}
where the summation is over all indecomposable partitions
$v=v_1\cup\cdots\cup v_p$ of the following two-way table
\begin{center}
\begin{tabular}{cc}
$r_1$&$r_1-j_1$\\
$r_2$&$r_2-j_2$\\
$r_1'$&$r_1'-j_1'$\\
$r_2'$&$r_2'-j_2'$.
\end{tabular}
\end{center}
Again by the absolute summability of $k$-th ($k=2,\cdots,8$)
cumulants, we get $H_{4n}(t,t')\le Cm_n^2 (t\vee t')^2$.
Therefore, $\E({J}_{33n}^2)\le Cm_n^3/n=o(m_n^2)$ and
${J}_{33n}=o_p(m_n)$. Thus the conclusion is established.

\qed

\begin{lemma}
\label{lem:J4n} Under the assumptions in Theorem~\ref{th:null},
the random variable $J_{4n}=4/n^2\sum_{t=6m_n+1}^{n}\sum_{r,r'=m_n+2}^{t-2l_n}\sum_{j,j'=1}^{m_n}k_{nj}^2k_{nj'}^2\E({D}_{j,t}{D}_{j',t})[\tilde{D}_{j,r}\tilde{D}_{j',r'}-\E(\tilde{D}_{j,r}\tilde{D}_{j',r'})]$ as defined in (\ref{eq:vnbar}) is $o_p(m_n)$.
\end{lemma}
\noindent Proof of Lemma~\ref{lem:J4n}: Write $J_{4n}=J_{41n}+J_{42n}$, where
\begin{eqnarray*}
J_{41n}&=&\frac{4}{n^2}\sum_{t=6m_n+1}^{n}\sum_{r,r'=m_n+2}^{t-2l_n}\sum_{j,j'=1,j\not=j'}^{m_n}k_{nj}^2k_{nj'}^2\E({D}_{j,t}{D}_{j',t})[\tilde{D}_{j,r}\tilde{D}_{j',r'}-\E(\tilde{D}_{j,r}\tilde{D}_{j',r'})],
\\
J_{42n}&=&\frac{4}{n^2}\sum_{t=6m_n+1}^{n}\sum_{r,r'=m_n+2}^{t-2l_n}\sum_{j=1}^{m_n}k_{nj}^4 \E({D}_{j,t}^2)[\tilde{D}_{j,r}\tilde{D}_{j,r'}-\E(\tilde{D}_{j,r}\tilde{D}_{j,r'})].
\end{eqnarray*}
Note that
\begin{eqnarray*}
&&\hspace{-0.5cm}\E(J_{41n}^2)=O(n^{-4})\sum_{t_1,t_2=6m_n+1}^{n}\sum_{r_1,r_1'=m_n+2}^{t_1-2l_n}\sum_{r_2,r_2'=m_n+2}^{t_2-2l_n}\sum_{j_1,j_1'=1,j_1\not=j_1'}^{m_n}\sum_{j_2,j_2'=1,j_2\not=j_2'}^{m_n}k_{nj_1}^2k_{nj_1'}^2\\
&&k_{nj_2}^2k_{nj_2'}^2\E({D}_{j_1,t_1}{D}_{j_1',t_1})\E({D}_{j_2,t_2}{D}_{j_2',t_2})\{\cov(\tilde{D}_{j_1,r_1},\tilde{D}_{j_2,r_2})\cov(\tilde{D}_{j_1',r_1'},\tilde{D}_{j_2',r_2'})\\
&&+\cov(\tilde{D}_{j_1,r_1},\tilde{D}_{j_2',r_2'})\cov(\tilde{D}_{j_1',r_1'},\tilde{D}_{j_2,r_2})+\cum(\tilde{D}_{j_1,r_1},\tilde{D}_{j_1',r_1'},\tilde{D}_{j_2,r_2},\tilde{D}_{j_2',r_2'})\}.
\end{eqnarray*}
By Lemma~\ref{lem:djk} and (\ref{eq:bound1}), the first two terms in the curly bracket above
contribute $O(m_n)$. Since
$\cum(\tilde{D}_{j_1,r_1},\tilde{D}_{j_1',r_1'},\tilde{D}_{j_2,r_2},\tilde{D}_{j_2',r_2'})$
vanishes when any two neighboring indices (say, $(r_1,r_1')$,
$(r_1',r_2)$ and $(r_2,r_2')$ if $r_1\ge r_1'\ge r_2\ge r_2'$) are
more than $l_n$ apart, the third term is $O(l_n^3/n)=o(m_n^2)$. So
$J_{41n}=o_p(m_n)$. Concerning $J_{42n}$, we have $J_{42n}=J_{421n}+J_{422n}$, where
\begin{eqnarray*}
J_{421n}&=&\frac{4}{n^2}\sum_{t=6m_n+1}^{n}\sum_{r=m_n+2}^{t-2l_n}\sum_{j=1}^{m_n}k_{nj}^4 \E({D}_{j,t}^2)[\tilde{D}_{j,r}^2-\E(\tilde{D}_{j,r}^2)],\\
J_{422n}&=&\frac{8}{n^2}\sum_{t=6m_n+1}^{n}\sum_{r=m_n+3}^{t-2l_n}\sum_{r'=m_n+2}^{r-1}\sum_{j=1}^{m_n}k_{nj}^4 \E({D}_{j,t}^2)[\tilde{D}_{j,r}\tilde{D}_{j,r'}-\E(\tilde{D}_{j,r}\tilde{D}_{j,r'})].
\end{eqnarray*}
Since $\tilde{D}_{j,r}^2$ is $l_n$-dependent, we can easily derive $\E(J_{421n}^2)=O(m_n^3/n)$, which implies $J_{421n}=o_p(m_n)$. Let
\[\tilde{J}_{422n}=\frac{8}{n^2}\sum_{t=6m_n+1}^{n}\sum_{r=m_n+3}^{t-2l_n}\sum_{r'=m_n+2}^{r-1}\sum_{j=1}^{m_n}k_{nj}^4 \E({D}_{j,t}^2){D}_{j,r}{D}_{j,r'}.\]
Then by (\ref{eq:bound1}), $J_{422n}-\tilde{J}_{422n}=o_p(1)$. Since for each $j$, $\{\sum_{r=m_n+3}^{t-2l_n}\sum_{r'=m_n+2}^{r-1}{D}_{j,r}{D}_{j,r'}\}$ form
martingale differences with respect to ${\cal F}_{t-2l_n}$, we get
\begin{eqnarray*}
\E(\tilde{J}_{422n}^2)&\le &Cn^{-4} m_n\sum_{j=1}^{m_n}k_{nj}^8\E \left(\sum_{t=6m_n+1}^{n}\sum_{r=m_n+3}^{t-2l_n}\sum_{r'=m_n+2}^{r-1}\E({D}_{j,t}^2){D}_{j,r}{D}_{j,r'}\right)^2\\
&\le&\frac{C m_n}{n^4}\sum_{j=1}^{m_n}\sum_{t=6m_n+1}^{n}\E\left(\sum_{r=m_n+3}^{t-2l_n}\sum_{r'=m_n+2}^{r-1}{D}_{j,r}{D}_{j,r'}\right)^2\\
&=&\frac{C m_n}{n^4}\sum_{j=1}^{m_n}\sum_{t=6m_n+1}^{n}\sum_{r=m_n+3}^{t-2l_n}\E\left[{D}_{j,r}^2\left(\sum_{r'=m_n+2}^{r-1}{D}_{j,r'}\right)^2\right],
\end{eqnarray*}
where we have applied the fact that for each $j$, $\{\sum_{r'=m_n+2}^{r-1}{D}_{j,r}{D}_{j,r'}\}$ is a sequence of martingale differences with respect to ${\cal F}_r$. By the Cauchy-Schwarz inequality and Burkholder's inequality,
\[\E\left[{D}_{j,r}^2\left(\sum_{r'=m_n+2}^{r-1}{D}_{j,r'}\right)^2\right]\le C\left\|\sum_{r'=m_n+2}^{r-1}{D}_{j,r'}\right\|_4^2\le C (r-m_n-2).\]
 Thus $\E(\tilde{J}_{422n}^2)\le Cm_n^2/n=o(m_n^2)$, in other words, $\tilde{J}_{422n}=o_p(m_n)$. The proof is complete.

\qed

\subsection{Appendix B}

Throughout the appendix B, we let $u_t(\theta)=\sum_{k=0}^{\infty}e_k(\theta)Y_{t-k}$ and
$\hat{u}_t=\sum_{k=0}^{t-1}e_k(\hat{\theta}_n)Y_{t-k}$,
$t=1,2,\cdots,n$. Write $\hat{u}_t=u_t+\lambda_{nt}$, where
$\lambda_{nt}=\lambda_{1t}+\lambda_{2nt}$,
$\lambda_{1t}=-\sum_{k=t}^{\infty}e_k(\theta_0)Y_{t-k}=\sum_{k=0}^{\infty}\psi_{k,t}u_{-k}$
and
$\lambda_{2nt}=\sum_{k=0}^{t-1}(e_k(\hat{\theta}_n)-e_k(\theta_0))Y_{t-k}$.
Denote by  $e_{k;m_1}(\theta)=\partial e_k(\theta)/\partial \theta_{m_1}$ and
$e_{k;(m_1,m_2)}(\theta)=\partial^2 e_k(\theta)/\partial
\theta_{m_1}\partial \theta_{m_2}$  for any $m_1,m_2\in
\{1,2,\cdots,p+q+1\}$ and assume they are the same as those expressions in Lemma~\ref{lem:bound}
without loss of generality.

\begin{lemma}
\label{lem:resi} Under the assumptions in
Theorem~\ref{th:farimaresidual}, we have (a).
$n\sum_{j=1}^{m_n}k_{nj}^2\hat{\rho}_{\hat{u}}^2(j)=n\sigma^{-4}\sum_{j=1}^{m_n}k_{nj}^2\hat{R}_{\hat{u}}^2(j)+o_p(m_n^{1/2})$.
and (b).
$n\sum_{j=1}^{m_n}k_{nj}^2(\hat{R}_{\hat{u}}^2(j)-\tilde{R}_{\hat{u}}^2(j))=o_p(m_n^{1/2})$,
where
$\tilde{R}_{\hat{u}}(j)=n^{-1}\sum_{t=|j|+1}^{n}\hat{u}_{t}\hat{u}_{t-|j|}$.
\end{lemma}
\noindent Proof of Lemma~\ref{lem:resi}:
To prove (a), it suffices to show that
\begin{eqnarray}
\label{eq:Ru0}
\hat{R}_{\hat{u}}(0)=n^{-1}\sum_{t=1}^{n}\hat{u}_t^2-\left(n^{-1}\sum_{t=1}^{n}\hat{u}_t\right)^2=\sigma^2+O_p(n^{-1/2}).
\end{eqnarray}
To this end, let $G_{1n}=n^{-1}\sum_{t=1}^{n}u_t\lambda_{1t}$,
$G_{2n}=n^{-1}\sum_{t=1}^{n}\lambda_{1t}^2$ and
$G_{3n}=n^{-1}\sum_{t=1}^{n}\lambda_{2nt}^2$. Since
$n^{-1}\sum_{t=1}^{n}u_t^2-\sigma^2=O_p(n^{-1/2})$, (\ref{eq:Ru0})
follows if we can show $G_{1n}=O_p(n^{-1/2})$,
$G_{2n}=O_p(n^{-1/2})$ and $G_{3n}=O_p(n^{-1})$.
Note that
\begin{eqnarray*}
\E(G_{1n}^2)&=&n^{-2}\sum_{t,t'=1}^{n}\sum_{k,k'=0}^{\infty}\psi_{k,t}\psi_{k',t'}\E(u_tu_{t'}u_{-k}u_{-k'})\\
 &=&n^{-2}\sum_{t=1}^{n}\sum_{k=0}^{\infty}\psi_{k,t}^2\sigma^4+n^{-2}\sum_{t,t'=1}^{n}\sum_{k,k'=0}^{\infty}\psi_{k,t}\psi_{k',t'}\cum(u_t,u_{t'},u_{-k},u_{-k'})\\
 &=&O(\log n/n^2+n^{-1})=O(n^{-1}),
\end{eqnarray*}
where we have applied the fact that
$\sum_{k=0}^{\infty}\psi_{k,t}^2=O(t^{-1})$ [cf. Robinson (2005)]
and the absolute summability of the 4-th cumulants.  Since
$\E(G_{2n})=O(\log n/n)$, $G_{2n}=O_p(n^{-1/2})$. To show $G_{3n}=o_p(n^{-1})$, we apply  the mean-value theorem and get
$e_k(\hat{\theta}_n)-e_k(\theta_0)=\sum_{m_1=1}^{p+q+1}e_{k;m_1}(\bar{\theta}_{kn})(\hat{\theta}_n^{(m_1)}-\theta_0^{(m_1)})$,
where
$\bar{\theta}_{kn}=\theta_0+\beta_k(\hat{\theta}_{n}-\theta_0)$
for some $\beta_k\in (0,1)$. Then
 \begin{eqnarray*}
&&n G_{3n}=\sum_{t=1}^{n}\sum_{k,k'=0}^{t-1}(e_k(\hat{\theta}_n)-e_k(\theta_0))(e_{k'}(\hat{\theta}_n)-e_{k'}(\theta_0))Y_{t-k}Y_{t-k'}\\
          &&=\sum_{t=1}^{n}\sum_{k,k'=0}^{t-1}\sum_{m_1,m_1'=1}^{p+q+1}(\hat{\theta}_n^{(m_1)}-\theta_0^{(m_1)})(\hat{\theta}_n^{(m_1')}-\theta_0^{(m_1')})e_{k;m_1}(\bar{\theta}_{kn})e_{k';m_1'}(\bar{\theta}_{k'n})Y_{t-k}Y_{t-k'}.
 \end{eqnarray*}
When $\hat{\theta}_n\in\Theta_{\delta}$, by Lemma~\ref{lem:bound},
for any $(m_1,m_1')\in\{1,\cdots,p+q+1\}^2$,
\[\sum_{t=1}^{n}\sum_{k,k'=0}^{t-1}|e_{k;m_1}(\bar{\theta}_{kn})||e_{k';m_1'}(\bar{\theta}_{k'n})|\E|Y_{t-k}Y_{t-k'}|=O(n).\]
Since $\hat{\theta}_n-\theta_0=O_p(n^{-1/2})$, we have
$P(\hat{\theta}_n\notin \Theta_{\delta})\rightarrow 0$.
Consequently $nG_{3n}=nG_{3n}{\bf 1}(\hat{\theta}_n\in \Theta_{\delta})+nG_{3n}{\bf 1}(\hat{\theta}_n\notin \Theta_{\delta})=O_p(1)$.
Therefore part (a) is proved.

As to part (b), write
$\tilde{R}_{\hat{u}}(j)-\hat{R}_{\hat{u}}(j)=-n^{-1}\bar{\hat{u}}\left(\sum_{t=1}^{n-j}\hat{u}_t+\sum_{t=j+1}^{n}\hat{u}_t\right)+
(1-j/n)\bar{\hat{u}}^2$,
where $\bar{\hat{u}}=n^{-1}\sum_{t=1}^{n}\hat{u}_t$. Following the
argument for part (a), it is straightforward to show that
$\bar{\hat{u}}=O_p(n^{-1/2})$ and
$\sum_{j=1}^{m_n}k_{nj}^2(\sum_{t=1}^{n-j}\hat{u}_t+\sum_{t=j+1}^{n}\hat{u}_t)^2=O_p(n m_n)$.
So
$n\sum_{j=1}^{m_n}k_{nj}^2(\hat{R}_{\hat{u}}(j)-\tilde{R}_{\hat{u}}(j))^2=o_p(1)$.
Applying the Cauchy-Schwarz inequality, part (b) follows. \qed

 \noindent {\it Proof of
Theorem~\ref{th:farimaresidual}:} By Lemma~\ref{lem:resi}, we only
need to show that
\[\frac{n\sum_{j=1}^{m_n}k_{nj}^2\tilde{R}_{\hat{u}}^2(j)-\sigma^4m_nC(K)}{(2\sigma^8 m_n D(K))^{1/2}}\rightarrow_{D}N(0,1).\]
 Note that
$\tilde{R}_{\hat{u}}^2(j)-\tilde{R}_u^2(j)=(\tilde{R}_{\hat{u}}(j)-\tilde{R}_u(j))^2+2\tilde{R}_u(j)(\tilde{R}_{\hat{u}}(j)-\tilde{R}_u(j))$.
By   Theorem~\ref{th:null}, it suffices to show
\begin{eqnarray*}
\label{eq:key}
n\sum_{j=1}^{m_n}k_{nj}^2(\tilde{R}_{\hat{u}}(j)-\tilde{R}_u(j))^2=o_p(1),
\end{eqnarray*}
since it implies
$n\sum_{j=1}^{m_n}k_{nj}^2\tilde{R}_u(j)(\tilde{R}_{\hat{u}}(j)-\tilde{R}_u(j))=o_p(m_n^{1/2})$
by the Cauchy-Schwarz inequality. To this end, we note that
\begin{eqnarray*}
&&n\sum_{j=1}^{m_n}k_{nj}^2(\tilde{R}_{\hat{u}}(j)-\tilde{R}_u(j))^2\le
\frac{C}{n}\sum_{j=1}^{m_n}\left\{\left(\sum_{t=j+1}^{n}\lambda_{1t}u_{t-j}\right)^2+
\left(\sum_{t=j+1}^{n}\lambda_{2nt}u_{t-j}\right)^2\right.\\
&&+\left(\sum_{t=j+1}^{n}u_t\lambda_{1(t-j)}\right)^2+\left(\sum_{t=j+1}^{n}u_t\lambda_{2n(t-j)}\right)^2+\left.\left(\sum_{t=j+1}^{n}\lambda_{nt}\lambda_{n(t-j)}\right)^2\right\}\\
&&\hspace{0.5cm}=: C (L_{1n}+L_{2n}+L_{3n}+L_{4n}+L_{5n}).
\end{eqnarray*}
We proceed to show that $L_{kn}=o_p(1), k=1,\cdots,5$. First,
\begin{eqnarray*}
\E(L_{1n})&=&
%n^{-1}\sum_{j=1}^{m_n}\E\left(\sum_{t=j+1}^{n}\lambda_{1t}u_{t-j}\right)^2=
n^{-1}\sum_{j=1}^{m_n}\E\left(\sum_{t=j+1}^{n}\sum_{k=0}^{\infty}\psi_{k,t}u_{-k}u_{t-j}\right)^2\\
&=&n^{-1}\sum_{j=1}^{m_n}\sum_{t,t'=j+1}^{n}\sum_{k=0}^{\infty}\sum_{k'=0}^{\infty}\psi_{k,t}\psi_{k',t'}\E(u_{-k}u_{-k'}u_{t-j}u_{t'-j})\\
&=&n^{-1}\sum_{j=1}^{m_n}\sum_{t,t'=j+1}^{n}\sum_{k=0}^{\infty}\sum_{k'=0}^{\infty}\psi_{k,t}\psi_{k',t'}\{\cov(u_{-k},u_{-k'})\cov(u_{t-j},u_{t'-j})\\
&&\hspace{1cm}+\cum(u_{-k},u_{-k'},u_{t-j},u_{t'-j})\},
\end{eqnarray*}
where the first term above is
$(\sigma^4/n)\sum_{j=1}^{m_n}\sum_{t=j+1}^{n}\sum_{k=0}^{\infty}\psi_{k,t}^2=O(m_n\log
n/n)$. Applying Proposition 2 in Wu and Shao (2004), we have
$|\cum(u_{-k},u_{-k'},u_{t-j},u_{t'-j})|\le C r^{t\vee t'-j+k\vee
k'}$ for some $r\in (0,1)$. So the second term in $\E(L_{1n})$ is
bounded by
\begin{eqnarray*}
%&&n^{-1}\sum_{j=1}^{m_n}\sum_{t,t'=j+1}^{n}\sum_{k=0}^{\infty}\sum_{k'=0}^{\infty}\psi_{k,t}\psi_{k',t'}|\cum(u_0,u_{k-k'},u_{k+t-j},u_{k+t'-j})|\\
%&&\hspace{0.5cm}\le
Cn^{-1}\sum_{j=1}^{m_n}\sum_{t,t'=j+1}^{n}\sum_{k=0}^{\infty}\sum_{k'=0}^{\infty}|\psi_{k,t}\psi_{k',t'}|r^{t\vee
t'-j+k\vee k'}=O(m_n/n).
\end{eqnarray*}
Following the same argument, we get $\E(L_{3n})=O(m_n/n)=o(1)$.

To show $L_{5n}=o_p(1)$, we note that
\begin{eqnarray}
\label{eq:L5n}
L_{5n}&\le&\frac{C}{n}\sum_{j=1}^{m_n}\left\{\left(\sum_{t=j+1}^{n}\lambda_{1t}\lambda_{1(t-j)}\right)^2+\left(\sum_{t=j+1}^{n}\lambda_{1t}\lambda_{2n(t-j)}\right)^2+\left(\sum_{t=j+1}^{n}\lambda_{2nt}\lambda_{1(t-j)}\right)^2\right.\nonumber\\
&&\left.+\left(\sum_{t=j+1}^{n}\lambda_{2nt}\lambda_{2n(t-j)}\right)^2\right\}=:C(L_{51n}+L_{52n}+L_{53n}+L_{54n}).
\end{eqnarray}
As to $L_{51n}$, we have
\begin{eqnarray*}
\E(L_{51n})&=&\frac{1}{n}\sum_{j=1}^{m_n}\sum_{t,t'=j+1}^{n}\E(\lambda_{1t}\lambda_{1t'}\lambda_{1(t-j)}\lambda_{1(t'-j)})\\
           &=&\frac{1}{n}\sum_{j=1}^{m_n}\sum_{t,t'=j+1}^{n}\sum_{k_1,k_2,k_3,k_4=0}^{\infty}\psi_{k_1,t}\psi_{k_2,t'}\psi_{k_3,t-j}\psi_{k_4,t'-j}\E(u_{-k_1}u_{-k_2}u_{-k_3}u_{-k_4})\\
           &=&\frac{1}{n}\sum_{j=1}^{m_n}\sum_{t,t'=j+1}^{n}\sum_{k_1,k_2,k_3,k_4=0}^{\infty}\psi_{k_1,t}\psi_{k_2,t'}\psi_{k_3,t-j}\psi_{k_4,t'-j}\{\cov(u_{-k_1},u_{-k_2})\\
           &&\cov(u_{-k_3},u_{-k_4})+\cov(u_{-k_1},u_{-k_3})\cov(u_{-k_2},u_{-k_4})\\
           &&+\cov(u_{-k_1},u_{-k_4})\cov(u_{-k_2},u_{-k_3})+\cum(u_{-k_1},u_{-k_2},u_{-k_3},u_{-k_4})\}.
\end{eqnarray*}
Since $\sum_{k=0}^{\infty}\psi_{k,t}^2\le Ct^{-1}$ [cf. Robinson
(2005)], the first three terms above are $O(m_n\log^2 n/n)$ under
the null hypothesis. By Proposition 2 in Wu and Shao (2004),
$|\cum(u_{-k_1},u_{-k_2},u_{-k_3},u_{-k_4})|\le
C r^{\max(k_1,k_2,k_3,k_4)-\min(k_1,k_2,k_3,k_4)}$ for some $r\in (0,1)$. Thus the
fourth term above is bounded by
\begin{eqnarray*}
&&\frac{C}{n}\sum_{j=1}^{m_n}\sum_{t,t'=j+1}^{n}\sum_{k_1\ge
k_2\ge k_3\ge
k_4=0}^{\infty}|\psi_{k_1,t}\psi_{k_2,t'}\psi_{k_3,t-j}\psi_{k_4,t'-j}| r^{k_1-k_4}\\
&&\hspace{0.5cm}\le
\frac{C}{n}\sum_{j=1}^{m_n}\sum_{t,t'=j+1}^{n}\sum_{h_1,h_3=0}^{\infty}\sum_{k_2=0}^{\infty}|\psi_{k_2+h_1,t}\psi_{k_2,t'}|
\sum_{k_4=0}^{\infty}|\psi_{k_4+h_3,t-j}\psi_{k_4,t'-j}| r^{h_1+h_3}\\
&&\hspace{0.5cm}\le
\frac{C}{n}\sum_{j=1}^{m_n}\sum_{t,t'=j+1}^{n}(tt'(t-j)(t'-j))^{-1/2}\sum_{h_1,h_3=0}^{\infty} r^{h_1+h_3}=o(1).
\end{eqnarray*}
Lemma~\ref{lem:L52n} asserts that $L_{52n}=o_p(1)$ and the same
argument leads to $L_{53n}=o_p(1)$. Following the line as in the derivation of $G_{3n}$ (see Lemma~\ref{lem:resi}),
we can derive
 $L_{54n}=O_p(m_n/n)=o_p(1)$. Thus $L_{5n}=o_p(1)$ and a
similar and simpler argument yields  $L_{kn}=o_p(1)$, $k=2,4$. We
omit the details. The conclusion is established.

\qed

\begin{lemma}
\label{lem:L52n} Under the assumptions in
Theorem~\ref{th:farimaresidual}, the random variable $L_{52n}=n^{-1}\sum_{j=1}^{m_n} \left(\sum_{t=j+1}^{n}\lambda_{1t}\lambda_{2n(t-j)}\right)^2$ as defined in (\ref{eq:L5n}) is $o_p(1)$.
%$L_{52n}$ defined in (\ref{eq:L5n}) is $o_p(1)$.
\end{lemma}
\noindent Proof of Lemma~\ref{lem:L52n}:  We apply  a Taylor's expansion for
each $k$ and obtain
\begin{eqnarray*}
\label{eq:taylor}
e_k(\hat{\theta}_n)-e_k(\theta_0)&=&\sum_{m_1=1}^{p+q+1}(\hat{\theta}_n^{(m_1)}-\theta_0^{(m_1)})e_{k;m_1}(\theta_0)\nonumber\\
&&\hspace{0.5cm}+\sum_{m_1,m_2=1}^{p+q+1}(\hat{\theta}_n^{(m_1)}-\theta_0^{(m_1)})
(\hat{\theta}_n^{(m_2)}-\theta_0^{(m_2)})e_{k;(m_1,m_2)}(\tilde{\theta}_{kn}),
\end{eqnarray*}
 where $\tilde{\theta}_{kn}=\theta_0+\alpha_k(\hat{\theta}_n-\theta_0)$
for some $\alpha_k\in (0,1)$. By Lemma~\ref{lem:bound},
% By Lemmas A1 and A2 in Francq and
%Zakoian (2000), we have that
$|e_{k;m_1}(\theta_0)|\le Ck^{-1-\epsilon}$ and
$\sup_{\theta\in\Theta_{\delta}}|e_{k;(m_1,m_2)}(\theta)|\le C
k^{-1-\epsilon}$ for some $\epsilon>0$. Denote by
$e_k(\theta_0)=e_k$ and $e_{k;m_1}(\theta_0)=e_{k;m_1}$. Since
$e_0(\theta)=1$, we have
\begin{eqnarray*}
L_{52n}&=&\frac{1}{n}\sum_{j=1}^{m_n}\sum_{t_1,t_2=j+1}^{n}\lambda_{1t_1}\lambda_{1t_2}\lambda_{2n(t_1-j)}\lambda_{2n(t_2-j)}\\
       &=&\frac{1}{n}\sum_{j=1}^{m_n}\sum_{t_1,t_2=j+1}^{n}\sum_{k_1,k_2=0}^{\infty}\sum_{k_3=1}^{t_1-j-1}\sum_{k_4=1}^{t_2-j-1}\psi_{k_1,t_1}\psi_{k_2,t_2}u_{-k_1}u_{-k_2}
       \\
       &&(e_{k_3}(\hat{\theta}_n)-e_{k_3}(\theta_0))(e_{k_4}(\hat{\theta}_n)-e_{k_4}(\theta_0))Y_{t_1-j-k_3}Y_{t_2-j-k_4}\\
       &=&\frac{1}{n}\sum_{j=1}^{m_n}\sum_{t_1,t_2=j+1}^{n}\sum_{k_1,k_2=0}^{\infty}\sum_{k_3=1}^{t_1-j-1}\sum_{k_4=1}^{t_2-j-1}\psi_{k_1,t_1}\psi_{k_2,t_2}u_{-k_1}u_{-k_2}Y_{t_1-j-k_3}Y_{t_2-j-k_4}\\
       &&\left(\sum_{m_1=1}^{p+q+1}(\hat{\theta}_n^{(m_1)}-\theta_0^{(m_1)})
       e_{k_3;m_1}+\sum_{m_1,m_2=1}^{p+q+1}(\hat{\theta}_n^{(m_1)}-\theta_0^{(m_1)})e_{k_3;(m_1,m_2)}(\tilde{\theta}_{k_3n})\right.\\
       &&\left.(\hat{\theta}_n^{(m_2)}-\theta_0^{(m_2)})\right)\left(\sum_{m_3=1}^{p+q+1}(\hat{\theta}_n^{(m_3)}-\theta_0^{(m_3)})
       e_{k_4;m_3}+\sum_{m_3,m_4=1}^{p+q+1}(\hat{\theta}_n^{(m_3)}-\theta_0^{(m_3)})\right.\\
       &&\left.e_{k_4;(m_3,m_4)}(\tilde{\theta}_{k_4n})(\hat{\theta}_n^{(m_4)}-\theta_0^{(m_4)})\right)=\sum_{h=1}^4 L_{52hn}.
\end{eqnarray*}
Write $Y_t=\sum_{k=0}^{\infty}a_k u_{t-k}$. To  show
$L_{521n}=o_p(1)$, it suffices to show that for any
$(m_1,m_3)\in\{1,\cdots,p+q+1\}^{2}$,
\begin{eqnarray*}
&&\tilde{L}_{521n}=\sum_{j=1}^{m_n}\sum_{t_1,t_2=j+1}^{n}\sum_{k_1,k_2=0}^{\infty}\psi_{k_1,t_1}\psi_{k_2,t_2}
\sum_{k_3=1}^{t_1-j-1}\sum_{k_4=1}^{t_2-j-1}
e_{k_3;m_1}e_{k_4;m_3}u_{-k_1}u_{-k_2}\\
&&\hspace{0.3cm}Y_{t_1-j-k_3}Y_{t_2-j-k_4}=\sum_{j=1}^{m_n}\sum_{t_1,t_2=j+1}^{n}\sum_{k_1,k_2=0}^{\infty}\sum_{k_3=1}^{t_1-j-1}\sum_{k_4=1}^{t_2-j-1}\sum_{h_1,h_2=0}^{\infty}\psi_{k_1,t_1}\psi_{k_2,t_2}\\
&&\hspace{0.3cm}a_{h_1}a_{h_2}e_{k_3;m_1}e_{k_4;m_3}u_{-k_1}u_{-k_2}u_{t_1-j-k_3-h_1}u_{t_2-j-k_4-h_2}=o_p(n^2).
\end{eqnarray*}
Note that
\begin{eqnarray*}
&&\E(\tilde{L}_{521n}^2)=\sum_{j,j'=1}^{m_n}\sum_{t_1,t_2=j+1}^{n}\sum_{t_1',t_2'=j'+1}^{n}\sum_{k_1,k_2,k_1',k_2'=0}^{\infty}\sum_{k_3=1}^{t_1-j-1}\sum_{k_4=1}^{t_2-j-1}\sum_{k_3'=1}^{t_1'-j'-1}\\
&&\hspace{0.3cm}\sum_{k_4'=1}^{t_2'-j'-1}\sum_{h_1,h_2,h_1',h_2'=0}^{\infty}\psi_{k_1,t_1}\psi_{k_2,t_2}
\psi_{k_1',t_1'}\psi_{k_2',t_2'}a_{h_1}a_{h_2}
a_{h_1'}a_{h_2'}e_{k_3;m_1}e_{k_4;m_3}
e_{k_3';m_1}e_{k_4';m_3}\\
&&\hspace{0.3cm}\E(u_{-k_1}u_{-k_2}u_{t_1-j-k_3-h_1}u_{t_2-j-k_4-h_2}u_{-k_1'}u_{-k_2'}u_{t_1'-j'-k_3'-h_1'}u_{t_2'-j'-k_4'-h_2'})\\
&&\hspace{0.3cm}\le C\sum_{j,j'=1}^{m_n}\sum_{t_1,t_2=j+1}^{n}\sum_{t_1',t_2'=j'+1}^{n}\sum_{k_1,k_2,k_1',k_2'=0}^{\infty}\sum_{k_3,k_4,k_3',k_4'=1}^{\infty}\sum_{h_1,h_2,h_1',h_2'=0}^{\infty}\\
&&\hspace{0.3cm}|\psi_{k_1,t_1}\psi_{k_2,t_2}||\psi_{k_1',t_1'}\psi_{k_2',t_2'}||a_{h_1}a_{h_2}||a_{h_1'}a_{h_2'}|
|k_3k_{3}'k_4k_4'|^{-1-\epsilon}\mbox{II},
\end{eqnarray*}
where
\begin{eqnarray*}
\mbox{II}&=&|\E(u_{-k_1}u_{-k_2}u_{t_1-j-k_3-h_1}u_{t_2-j-k_4-h_2}u_{-k_1'}u_{-k_2'}u_{t_1'-j'-k_3'-h_1'}u_{t_2'-j'-k_4'-h_2'})|\\
&=&\left|\sum_{g}\cum(u_{i_j},i_j\in g_1)\cdots
\cum(u_{i_j},i_j\in g_p)\right|.
\end{eqnarray*}
 In the above equation, $\Sigma_g$ is over all partitions $g=\{g_1\cup
\cdots\cup g_p\}$ of the index set
$\{-k_1,t_1-j-k_3-h_1,-k_1',t_1'-j'-k_3'-h_1',-k_2,t_2-j-k_4-h_2,-k_2',t_2'-j'-k_4'-h_2'\}$.
Since $\E(u_t)=0$, only partitions $g$ with $\#g_i>1$ for all $i$
contribute. We shall divide all contributing partitions into the
following several types and treat them one by one.
\begin{enumerate}
\item $\#g_1=\#g_2=\#g_3=\#g_4=2$. One such term is
\begin{eqnarray*}
\label{eq:2222}
&&\cov(u_{-k_1},u_{t_1-j-k_3-h_1})\cov(u_{-k_1'},u_{t_1'-j'-k_3'-h_1'})\cov(u_{-k_2},u_{t_2-j-k_4-h_2})\\
&&\times \cov(u_{-k_2'},u_{t_2'-j'-k_4'-h_2'}),
\end{eqnarray*}
which is nonzero when $-k_1=t_1-j-k_3-h_1$,
$-k_1'=t_1'-j'-k_3'-h_1'$, $-k_2=t_2-j-k_4-h_2$ and
$-k_2'=t_2'-j'-k_4'-h_2'$. Define $a_h=0$ if $h<0$. Then for any
fixed $g\in\Z$, $\sum_{h=0}^{\infty}|a_ha_{h+g}|\le
\sum_{h=0}^{\infty}a_h^2:=S_a<\infty$. For any fixed
$t_1,t_1',t_2,t_2',j,j',k_3,k_4,k_3',k_4'$, by the Cauchy-Schwarz
inequality,
\begin{eqnarray*}
&&\sum_{k_1,k_2,k_1',k_2'=0}^{\infty}|\psi_{k_1,t_1}\psi_{k_2,t_2}||\psi_{k_1',t_1'}\psi_{k_2',t_2'}||a_{k_1+t_1-j-k_3}a_{k_2+t_2-j-k_4}|\\
&&|a_{k_1'+t_1'-j'-k_3'}a_{k_2'+t_2'-j'-k_4'}|\le
\left(\sum_{k_1,k_2,k_1',k_2'=0}^{\infty}\psi_{k_1,t_1}^2\psi_{k_2,t_2}^2\psi_{k_1',t_1'}^2\psi_{k_2',t_2'}^2\right)^{1/2}
S_a^{2}\\
&&\hspace{1cm}=O((t_1t_2t_1't_2')^{-1/2}).
\end{eqnarray*}
So this term is $O(m_n^2n^2)=o(n^4)$. Similarly,  all
non-vanishing terms involve four restrictions on the indices
$k_1,k_2,k_1',k_2',h_1,h_2,h_1',h_2'$ once we fix $t_1$ , $t_1'$,
$t_2$, $t_2'$, $j$, $j'$, $k_3$, $k_4$, $k_3'$, $k_4'$. The
contribution from these terms are of order $o(n^4)$.

\item $\#g_1=\#g_2=3,\#g_3=2$. A typical term is
\begin{eqnarray*}
&&\cum(u_{-k_1},u_{t_1-j-k_3-h_1},u_{-k_1'})\cum(u_{t_1'-j'-k_3'-h_1'},u_{-k_2},u_{t_2-j-k_4-h_2})\\
&&\times \cov(u_{-k_2'},u_{t_2'-j'-k_4'-h_2'}).
\end{eqnarray*}
So for any fixed $t_1,t_1',t_2,t_2',j,j',k_3,k_4,k_3',k_4'$,
\begin{eqnarray}
\label{eq:ii}
&&\hspace{1cm}\sum_{k_1,k_1',h_1=0}^{\infty}|\psi_{k_1,t_1}\psi_{k_1',t_1'}a_{h_1}||\cum(u_{-k_1},u_{t_1-j-k_3-h_1},u_{-k_1'})|\\
&\le&C\sum_{k_1,k_1',h_1=0}^{\infty}|\psi_{k_1,t_1}\psi_{k_1',t_1'}a_{h_1}| r^{\max(-k_1,t_1-j-k_3-h_1,-k_1')-\min(-k_1,t_1-j-k_3-h_1,-k_1')}.\nonumber
\end{eqnarray}
Consider the case $-k_1'\ge -k_1\ge t_1-j-k_3-h_1$. Then the
corresponding term above is
\[C\sum_{s_1,s_2,k_1=0}^{\infty}|\psi_{k_1,t_1}\psi_{k_1-s_1,t_1'}a_{s_2+k_1+t_1-j-k_3}| r^{s_1+s_2}=O((t_1t_1')^{-1/2}), \]
where we have applied the Cauchy-Schwarz inequality and the fact
that $\sum_{k=0}^{\infty}\psi_{k,t}^2=O(t^{-1})$. Other cases can
be treated in a similar fashion. So (\ref{eq:ii}) is
$O((t_1t_1')^{-1/2})$.  Similarly, we can show that
\[\sum_{h_1',h_2,k_2=0}^{\infty}|\psi_{k_2,t_2}a_{h_1'}a_{h_2}\cum(u_{t_1'-j'-k_3'-h_1'},u_{-k_2},u_{t_2-j-k_4-h_2})|=O(t_2^{-1/2})\]
and
\[\sum_{k_2',h_2'=0}^{\infty}|a_{h_2'}\psi_{k_2',t_2'}\cov(u_{-k_2'},u_{t_2'-j'-k_4'-h_2'})|=O((t_2')^{-1/2}).\]
 Thus these terms contribute $O(m_n^2n^2)=o(n^4)$.

\item $\#g_1=\#g_2=4$; $\#g_1=4,\#g_2=\#g_3=2$; $\#g_1=5,\#g_2=3$;
$\#g_1=6,\#g_2=2$ and $\#g_1=8$. Following a similar argument as
the second case, it is not hard to see that the contribution of
all these terms  are $o(n^4)$.
%or smaller due to the summability of the joint cumulants.
\end{enumerate}

So $L_{521n}=o_p(1)$. Under the assumption that $u_t$ is GMC$(8)$, it is not hard to show that $\E(Y_t^4)<\infty$, and $\sup_{t\in\N}\E\lambda_{1t}^4<\infty$; compare the derivation of
$\E(L_{51n})$ in the proof of Theorem~\ref{th:farimaresidual}. Together with  Lemma~\ref{lem:bound}, we have
$\E|L_{522n}|{\bf
1}(\hat{\theta}_n\in\Theta_{\delta})=O(m_n/n^{1/2})=o(1)$,
so $L_{522n}=o_p(1)$.  Similarly we derive  $L_{52kn}=o_p(1)$,
$k=3,4$. Now the proof is complete.

\qed

The following lemma is an extension of Lemma A.1 of Francq and
Zako\"1an (2000) to the FARIMA model.
\begin{lemma}
\label{lem:bound} For any $\theta\in\Theta_{\delta}$ and any $(m_1,m_2)\in
\{1,\cdots,p+q+1\}^2$, there exist absolutely summable sequences
$(e_{k}(\theta))_{k\ge 0}$, $(e_{k;m_1}(\theta))_{k\ge 1}$ and
$(e_{k;(m_1,m_2)}(\theta))_{k\ge 1}$
 such that almost surely
\begin{eqnarray*}
u_t(\theta)=\sum_{k=0}^{\infty}e_k(\theta)Y_{t-k}, ~\frac{\partial
u_{t}(\theta)}{\partial\theta_{m_1}}=\sum_{k=1}^{\infty}e_{k;m_1}(\theta)Y_{t-k}
\end{eqnarray*}
and
\[\frac{\partial^2u_t(\theta)}{\partial\theta_{m_1}\partial\theta_{m_2}}=\sum_{k=1}^{\infty}e_{k;(m_1,m_2)}(\theta)Y_{t-k}\]
Further,
there exists an $\epsilon>0$, such that
\[\sup_{\theta\in\Theta_{\delta}}|e_k(\theta)|=O(k^{-1-\epsilon}), ~\sup_{\theta\in\Theta_{\delta}}|e_{k;m_1}(\theta)|=O(k^{-1-\epsilon}),~\mbox{and}\]
\[\sup_{\theta\in\Theta_{\delta}}|e_{k;(m_1,m_2)}(\theta)|=O(k^{-1-\epsilon}).\]
\end{lemma}

\noindent Proof of Lemma~\ref{lem:bound}: Letting
$X_t=(1-B)^dY_t$, then $\phi_{\Lambda}(B)X_t=\psi_{\Lambda}(B)u_t$.
 By Lemma A.1 in Francq and Zako\"1an (2000), there exist
sequences $(c_{k}(\Lambda))_{k\ge 0}$, $(c_{k;m_1}(\Lambda))_{k\ge
1}$ and $(c_{k;(m_1,m_2)}(\Lambda))_{k\ge 1}$ such that
\[u_t(\Lambda)=\sum_{j=0}^{\infty}c_{j}(\Lambda)X_{t-j}, ~\partial
u_t(\Lambda)/\partial\Lambda_{m_1}=\sum_{j=1}^{\infty}c_{j;m_1}(\Lambda)X_{t-j}\]
and
\[\partial^2
u_t(\Lambda)/\partial\Lambda_{m_1}\partial\Lambda_{m_2}=\sum_{j=1}^{\infty}c_{j;(m_1,m_2)}(\Lambda)X_{t-j}.\]
Further, there exists a $r\in [0,1)$, such that
\[\sup_{\Lambda\in\Omega_{\delta}}|c_j(\Lambda)|=O(r^j),~\sup_{\Lambda\in\Omega_{\delta}}|c_{j;m_1}(\Lambda)|=O(r^j),~\sup_{\Lambda\in\Omega_{\delta}}|c_{j;(m_1,m_2)}(\Lambda)|=O(r^j).\]
Note
that $X_t=\sum_{s=0}^{\infty}\phi_s(d)Y_{t-s}$, where
$\phi_s(d)=\Gamma(s-d)/\{\Gamma(-d)\Gamma(s+1)\}$. Therefore,
we get $u_{t}(\theta)=\sum_{k=0}^{\infty}e_k(\theta)Y_{t-k}$,
where $e_{k}(\theta)=\sum_{j=0}^{k}c_{j}(\Lambda)\phi_{k-j}(d)$.
The conclusion follows from the definition of $\Theta_{\delta}$
and the fact that $d_0\in (0,1/2)$.

\qed

\end{document}